\documentclass[11pt]{article}
\usepackage{amsfonts}
\usepackage{amsmath}
\usepackage{epsfig}
\usepackage[compress]{cite}
\usepackage{graphics}
\usepackage{algorithm}
\usepackage{algorithmic}
  
\begin{document}

\begin{center}
\Large \bf{Implicit particle filtering for models with partial noise, and an application to geomagnetic data assimilation} 
\end{center}

\begin{center}
 Matthias Morzfeld$^1$ and Alexandre J. Chorin$^{1,2}$\\
\vspace{2mm}
$^1$Lawrence Berkeley National Laboratory, Berkeley, CA\\
and\\
$^2$Department of Mathematics\\
University of California, Berkeley, CA\\
\end{center}
\vspace{1mm}

\begin{center}
\bf{Abstract}
\end{center}
Implicit particle filtering is a sequential Monte Carlo method for data assimilation, designed to keep the number of particles manageable by focussing attention on regions of large probability. These regions are found by minimizing, for each particle, a scalar function $F$ of the state variables. Some previous implementations of the implicit filter rely on finding the Hessians of these functions. The calculation of the Hessians can be cumbersome if the state dimension is large or if the underlying physics are such that derivatives of $F$ are difficult to calculate. This is the case in many geophysical applications, in particular for models with partial noise, i.e. with a singular state covariance matrix. Examples of models with partial noise include stochastic partial differential equations driven by spatially smooth noise processes and models for which uncertain dynamic equations are supplemented by conservation laws with zero uncertainty.
We make the implicit particle filter applicable to such situations by combining gradient descent minimization with random maps and show that the filter is efficient, accurate and reliable because it operates in a subspace whose dimension is smaller than the state dimension. As an example, we assimilate data for a system of nonlinear partial differential equations that appears in models of geomagnetism.
\vspace{5mm}
\noindent



\section{Introduction}
\label{sec:Introduction}
The task in data assimilation is to use available data to update the forecast of a numerical model. The numerical model is typically given by a discretization of a stochastic differential equation (SDE)
\begin{equation}
	\label{eq:DSDE}
	x^{n+1}=R(x^n,t^n)+G(x^n,t^n)\Delta W^{n+1},
\end{equation}
where $x$ is an $m$-dimensional vector, called the state, $t^n$, $n=0,1,2,\dots$, is a sequence of times, $R$ is an $m$-dimensional vector function, $G$ is an $m\times m$ matrix and $\Delta W$ is an $m$-dimensional vector, whose elements are independent standard normal variates. The random vectors $G(x^n,t^n)\Delta W^{n+1}$ represent the uncertainty in the system, however even for $G=0$ the state $x^n$ may be random for any $n$ because the initial state $x^0$ can be random. The data 
\begin{equation}
	   \label{eq:Obs}
	z^{l}=h(x^{q(l)},t^{q(l)})+Q(x^{q(l)},t^{q(l)}) V^{l},
\end{equation}
are collected at times $t^{q(l)}$, $l=1,2,\dots$; for simplicity, we assume that the data are collected at a subset of the model steps, i.e. $q(l)=rl$, with $r\geq 1$ being a constant. In the above equation, $z$ is a $k$-dimensional vector ($k\leq m$), $h$ is a $k$-dimensional vector function, $V$ is a $k$-dimensional vector whose components are independent standard normal variates, and $Q$ is a $k \times k$ matrix. Throughout this paper, we will write $x^{1:n}$ for the sequence of vectors $\left\{x^0,\dots,x^n\right\}$.

Data assimilation is necessary in many areas of science and engineering and is essential in geophysics, for example in oceanography, meteorology, geomagnetism or atmospheric chemistry (see e.g. the reviews \cite{Miller1994,Ide1997,Miller1999,vanLeeuwen2009,Bocquet2010, Fournier2010}). What makes the assimilation of data in geophysical applications difficult is the complicated underlying physics, which lead to a large state dimension $m$ and a nonlinear function $R$ in equation~(\ref{eq:DSDE}).

If the model (\ref{eq:DSDE}) as well as $h$ in (\ref{eq:Obs}) are linear and if, in addition, the initial state $x^0$ is Gaussian, then the probability density function (pdf) of the state $x^n$ is Gaussian for any $n$ and can be characterized in full by its mean and covariance. The Kalman filter (KF) sequentially computes the mean of the model (\ref{eq:DSDE}), conditioned on the observations and, thus, provides the best linear unbiased estimate of the state \cite{Kalman1961}. The ensemble Kalman filter (EnKF) is a Monte Carlo approximation of the Kalman filter and can be obtained by replacing the state covariance matrix by the sample covariance matrix in the Kalman formalism. The state covariance is the covariance matrix of the pdf of the current state conditioned on the previous state which we calculate from the model (\ref{eq:DSDE}) to be:
\begin{equation}
	p(x^{n+1}\mid x^n)\sim \mathcal{N}(R(x^n,t^n),G(x^n,t^n)G(x^n,t^n)^T),
\end{equation}
where $\mathcal{N}(\mu,\Sigma)$ denotes a Gaussian with mean $\mu$ and covariance matrix $\Sigma$. To streamline the notation we write for the state covariance
\begin{equation}
	\label{eq:StateCovariance}
	\Sigma_x^n = G(x^n,t^n)G(x^n,t^n)^T,
\end{equation}
where $T$ denotes a transpose. In the EnKF, the sample covariance matrix is computed from an ``ensemble,'' by running the model (\ref{eq:DSDE}) for different realizations of the noise process $\Delta W$. The Monte Carlo approach avoids the computationally expensive step of updating the state covariance in the Kalman formalism. Both KF and EnKF have extensions to nonlinear, non-Gaussian models, however they rely on linearity and Gaussianity approximations \cite{Julier1997}. 

Variational methods \cite{Zupanski1997,Tremolet2006,Talagrand1997,Courtier1997,Courtier1994,Bennet1993,Talagrand1987} aim at assimilating the observations within a given time window by computing the state trajectory of maximum probability. The trajectory is computed by minimizing a suitable cost function which is, up to a normalization constant, the logarithm of the pdf of the state trajectory $=x^{0:q(l)}$ given the set of observations $z^{1:l}$, $p(x^{0:q(l)}\mid z^{1:l})$. In particular, 3D-Var methods assume a static model \cite{Talagrand1997}. Strong constraint 4D-Var determines an optimal initial state given a ``perfect'' dynamic model, i.e. $G=0$, and a Gaussian initial uncertainty, i.e. $x^0\sim\mathcal{N}(\mu^0,\Sigma^0)$ \cite{Talagrand1997,Courtier1997,Courtier1994,Talagrand1987}. Uncertain models with $G\neq 0$ are tackled with a weak constraint 4D-Var approach \cite{Zupanski1997,Tremolet2006,Bennet1993}. Many implementations of variational methods compute the gradient of the cost function from tangent linear adjoint equations and rely on linear approximations.

For the reminder of this paper, we focus on sequential Monte Carlo (SMC) methods for data assimilation, called particle filters \cite{Doucet2001,Weare2009,DelMoral1998,VanLeeuwen2010,DelMoral2004,Arulampalam2002,Doucet2000,Chorin2009,Chorin2010,Gordon1993,Morzfeld2011}. Particle filters do not rely upon linearity or Gaussianity assumptions and approximate the pdf of the state given the observations, $p(x^{0:q(l)}\mid z^{1:l})$, by SMC. The state estimate is a statistic (e.g. the mean, median, mode etc.) of this pdf. Most particle filters rely on the recursive relation
 \begin{equation}
	\label{eq:Recursive}
	p(x^{0:q(l+1)}\mid z^{1:l+1})\propto p(x^{0:q(l)}\mid z^{1:l})p(z^{l+1}\mid x^{q(l+1)})p(x^{q(l)+1:q(l+1)}\mid x^{q(l)}).
\end{equation}
In the above equation $p(x^{0:q(l+1)}\mid z^{1:l+1})$ is the pdf of the state trajectory up to time $t^{q(l+1)}$ given all available observations up to time $t^{q(l+1)}$ and is called the target density; $p(z^{l+1}\mid x^{q(l+1)})$ is the probability density of the current observation given the current state and can be obtained from (\ref{eq:Obs}):
\begin{equation}
	p(z^{l+1}\mid x^{q(l+1)})\sim\mathcal{N}(h(x^{q(l)},t^{q(l)}),\Sigma_z^{n}),
\end{equation}	
with
\begin{equation}
\Sigma_z^{n} = Q(x^n,t^n)Q(x^n,t^n)^T.
\end{equation}
The pdf $p(x^{q(l)+1:q(l+1)}\mid x^{q(l)})$ is the density of the state trajectory from the previous assimilation step to the current observation, conditioned on the state at the previous assimilation step, and is determined by the model (\ref{eq:DSDE}).

A standard version of the sampling-importance-resampling (SIR) particle filter (also called bootstrap filter, see e.g. \cite{Doucet2001}) generates, at each step, samples from $p(x^{q(l)+1:q(l+1)}\mid x^{q(l)})$ (the prior density) by running the model. These samples (particles) are weighted by the observations with weights $w\propto p(z^{l+1}\mid x^{q(l+1)})$, to yield a posterior density that approximates the target density $p(x^{0:q(l)}\mid z^{1:l})$. One then removes particles with a small weight by ``resampling'' (see e.g. \cite{Arulampalam2002} for resampling algorithms) and repeats the procedure when the next observation becomes available. This SIR filter is straightforward to implement, the catch is that many particles have small weights because the particles are generated without using information from the data. If many particles have a small weight, the approximation of the target density is poor and the number of particles required for a good approximation of the target density can grow catastrophically with the dimension of the state \cite{Snyder2008,Bickel2008}. Various methods, e.g. different prior densities and weighting schemes, have been invented to ameliorate this problem (see e.g. \cite{Doucet2001, VanLeeuwen2010,vanLeeuwen2009,Weare2009}).

The basic idea of implicit particle filters \cite{Chorin2009,Chorin2010,Morzfeld2011} is to use the available observations to find regions of high probability in the target density and look for samples within this region. This implicit sampling strategy generates a thin particle beam within the high probability domain and, thus, keeps the number of particles required manageable, even if the state dimension is large. The focussing of particles is achieved by setting up an underdetermined algebraic equation that depends on the model (\ref{eq:DSDE}) as well as on the data (\ref{eq:Obs}), and whose solution generates a high probability sample of the target density. We review the implicit filter in the next section, and it will become evident that the construction assumes that the state covariance $\Sigma_x^n$ in (\ref{eq:StateCovariance}) is nonsingular. This condition is often not satisfied. If, for example, one wants to assimilate data into a stochastic partial differential equation (SPDE) driven by spatially smooth noise, then the continuous-time noise process can be represented by a series with rapidly decaying coefficients, leading to a non-singular or ill-conditioned state covariance $\Sigma_x^n$ in discrete time and space (see Sections \ref{sec:Heat} and \ref{sec:Geomagnetism}, as well as \cite{LordRougemont,Chueshov2000,Jentzen2009}). A second important class of models with partial noise are uncertain dynamic equations supplemented by conservation laws (e.g. conservation of mass) with zero uncertainty. Such models often appear in data assimilation for fluid dynamics problems \cite{Kurapov2007}.

The purpose of the present paper is two-fold. First, in Section \ref{sec:ImplicitSampling}, we present a new implementation of the implicit particle filter. Most previous implementations of the implicit filter \cite{Chorin2010,Morzfeld2011} rely in one way or another on finding the Hessians of scalar functions in $rm$ variables. For systems with very large state vectors and considerable gaps between observations, memory constraints may forbid a computation of these Hessians. Our new implementation combines gradient descent minimization with random maps \cite{Morzfeld2011} to avoid the calculation of Hessians, and thus reduces the memory requirements. 

The second objective is to consider models with a singular or ill-conditioned state covariance $\Sigma_x^n$ where previous implementations of the implicit filter, as described in \cite{Chorin2009, Chorin2010,Morzfeld2011}, are not applicable. In Section \ref{sec:ImplicitSamplingPartialNoise}, we make the implicit filter applicable to models with partial noise and show that our approach is then particularly efficient, because the filter operates in a space whose dimension is determined by the rank of $\Sigma_x^n$, rather than by the model dimension. We compare the new implicit filter to SIR, EnKF and variational methods, in particular with respect to how information is propagated from observed variables to unobserved ones. 

In Section \ref{sec:Geomagnetism}, we illustrate the theory with an application in geomagnetism and consider two coupled nonlinear SPDE's with partial noise. We observe that the implicit filter gives good results with very few (4-10) particles, while EnKF and SIR require hundreds to thousands of particles for similar accuracy.

\section{Implicit sampling with random maps}\label{sec:ImplicitSampling}
We first follow \cite{Morzfeld2011} closely to review implicit sampling with random maps. Suppose we are given a collection of $M$ particles $X_j^{q(l)}$, $j=1,2,\dots,M$, whose empirical distribution approximates the target density at time $t^{q(l)}$, where $q(l)=rl$, and suppose that an observation $z^{l+1}$ is available after $r$ steps at time $t^{q(l+1)}=t^{r(l+1)}$. From (\ref{eq:Recursive}) we find, by repeatedly using Bayes' theorem, that, for each particle,
 \begin{eqnarray}
	p(X_j^{0:q(l+1)}\mid z^{1:l+1})&\propto& p(X_j^{0:q(l)}\mid z^{1:l}) p(z^{l+1}\mid X_j^{q(l+1)})\nonumber\\
									   & \times& p(X_j^{q(l+1)}\mid X_j^{q(l+1)-1}) p(X_j^{q(l+1)-1}\mid X_j^{q(l+1)-2}) \nonumber \\
									   & & \vdots \nonumber \\
					&\times &	 p(X_j^{q(l)+1}\mid X_j^{q(l)}). 
\end{eqnarray}
Implicit sampling is a recipe for computing high-probability samples from the above pdf. To draw a sample we define, for each particle, a function $F_j$ by
\begin{eqnarray}
	\label{eq:FDef}
	\exp(-F(X_j))&=&  p(X_j^{q(l+1)}\mid X_j^{q(l+1)-1}) \cdots  
			   p(X_j^{q(l)+1}\mid X_j^{q(l)})  \nonumber\\
			   &\times& p(z^{l+1}|X_j^{q(l+1)}).   
\end{eqnarray}
where $X_j$ is shorthand for the state trajectory $X_j^{q(l)+1:q(l+1)}$. Specifically, we have 
\begin{eqnarray}
	   \label{eq:FFullNoise}
	F_j(X_j) &=&\frac{1}{2} \left(X_j^{q(l)+1}-R_j^{q(l)}\right)^T \left(\Sigma_{x,j}^{q(l)}\right)^{-1}\left(X_j^{q(l)+1}-R_j^{q(l)}\right) \nonumber\\
			 & +&\frac{1}{2} \left(X_j^{q(l)+2}-R_j^{q(l)+1}\right)^T \left(\Sigma_{x,j}^{q(l)+1}\right)^{-1}\left(X_j^{q(l)+2}-R_j^{q(l)+1}\right) \nonumber\\
			 &  & \vdots \nonumber \\
			&+ &\frac{1}{2}  \left(X_j^{q(l+1)}-R_j^{q(l+1)-1}\right)^T \left(\Sigma_{x,j}^{q(l+1)-1}\right)^{-1}\left(X_j^{q(l+1)}-R_j^{q(l+1)-1}\right) \nonumber\\
			&+ &\frac{1}{2}  \left(h\left(X_j^{q(l+1)}\right)-z^{l+1}\right)^T \left(\Sigma_{z,j}^{l+1}\right)^{-1} \left(h\left(X_j^{q(l+1)}\right)-z^{l+1}\right) \nonumber \\
			&+&Z_j,
\end{eqnarray}
where $R_j^n$ is shorthand notation for $R(X_j^n,t^n)$ and where $Z_j$ is a positive number that can be computed from the normalization constants of the various pdf's in the definition of $F_j$ in (\ref{eq:FDef}).
With this $F_j$, we solve the algebraic equation
\begin{equation}
	\label{eq:AlgEq}
	F(X_j)-\phi_j = \frac{1}{2}\xi_j^T\xi_j,
\end{equation}
where $\xi_j$ is a realization of the $rm-$dimensional reference variable $\xi\sim\mathcal{N}(0,I)$, and where 
\begin{equation}
	\label{eq:Phi}
	\phi_j = \min F_j.
\end{equation}
The choice of a Gaussian reference variables does not imply linearity or Gaussianity assumptions and other choices are certainly possible. We find solutions of (\ref{eq:AlgEq}) by using the random map
\begin{equation}
        \label{eq:RandomMap}
	X_j = \mu_j + \lambda_j L_j \eta_j,
\end{equation}
where $\lambda_j$ is a scalar, $\mu_j$ is an $rm$-dimensional column vector which represents the location of the minimum of $F_j$, i.e. $\mu_j=\mbox{argmin} F_j$, $L_j$ is a deterministic $rm \times rm$ matrix we can choose, and $\eta_j=\xi_j/\sqrt{\xi_j^T\xi_j}$, is uniformly distributed on the unit $rm$-sphere. Upon substitution of (\ref{eq:RandomMap}) into (\ref{eq:AlgEq}), we can find a solution of (\ref{eq:AlgEq}) by solving a single algebraic equation in the variable $\lambda_j$. The weight of the particle can be shown to be (see \cite{Morzfeld2011}, Section 3)
\begin{equation}
        \label{eq:Weight}
	w_j^{q(l+1)} \propto w_j^{q(l)}\exp(-\phi_j) \left|\det L_j\right|  \mbox{ } \rho_j^{1-rm/2} \mbox{ } \left|\lambda_j^{rm-1} \frac{\partial \lambda_j}{\partial \rho_j}\right|,
\end{equation}
where $\rho_j = \xi_j^T\xi_j$ and $\det L_j$ denotes the determinant of the matrix $L_j$ (see \cite{Morzfeld2011} for details of the calculation). An expression for the scalar derivative $\partial\lambda_j/\partial\rho_j$ can be obtained by implicit differentiation of (\ref{eq:AlgEq}):
\begin{equation}
	\frac{\partial \lambda_j}{\partial \rho_j}= \frac{\rho}{2\left(\nabla F_j\right) L_j^T\eta_j},
\end{equation}
where $\nabla  F_j $ denotes the gradient of $F_j$ (an $rm$-dimensional row vector).

The weights are normalized so that their sum equals one. The weighted positions $X_j$ of the particles approximate the target pdf. We compute the mean of $X_j$ with weights $w_j$ as the state estimate, and then proceed to assimilate the next observation. The method just described makes use of only one set of observations per assimilation step, however an extension to multiple observation sets per assimilation step (smoothing) is straightforward.

\subsection{Implementation of an implicit particle filter with gradient descent minimization and random maps}
An algorithm for data assimilation with implicit sampling and random maps was presented in \cite{Morzfeld2011}. This algorithm relies on the calculation of the Hessians of the $F_j$'s and the Hessians are used for minimizing the $F_j$'s with Newton's method and for setting up the random map. The calculation of the Hessians however may not be easy in some applications because of a very large state dimension, or because the second derivatives are hard to calculate, as is the case for models with partial noise (see Section \ref{sec:ImplicitSamplingPartialNoise}). To avoid the calculation of Hessians, we propose to use a gradient descent algorithm with line-search to minimize the $F_j$'s (see e.g. \cite{Nocedal}), along with simple random maps.
Of course other minimization techniques, in particular quasi-Newton methods (see e.g. \cite{Nocedal,Fletcher}, can also be applied here and perhaps speed up the minimization. However, we decided to use gradient descent with line search to keep the minimization as simple as possible.

For simplicity, we assume that $G$ and $Q$ in (\ref{eq:DSDE})-(\ref{eq:Obs}) are constant matrices and calculate the gradient of $F_j$ from (\ref{eq:FFullNoise}):
\begin{equation}
	\label{eq:GradFullNoise}
	\nabla F = \left( 
					\frac{\partial F}{\partial X^{q(l)+1}},
					\frac{\partial F}{\partial X^{q(l)+2}},
					\dots, 
					\frac{\partial F}{\partial X^{q(l+1)-1}},
					\frac{\partial F}{\partial X^{q(l+1)}}
					\right),
\end{equation}
with
\begin{eqnarray}
	\left(\frac{\partial F}{\partial X^{k}}\right)^T &=& \Sigma_x^{-1} \left(X^{k}-R^{k-1} \right) \nonumber \\
								         & -&( \frac{\partial R}{\partial x} \mid _{x = X^{k}})^T
								         		              \Sigma_x^{-1}\left(X^{k+1}-R^{k}\right),
\end{eqnarray}
for $k=q(l)+1,q(l)+2,\dots,q(l+1)-1$, where $R^n$ is shorthand for $R(X^n,t^n)$, and where
\begin{eqnarray}
		\label{eq:LastGradFullNoise}
	\left(\frac{\partial F}{\partial X^{q(l+1)}}\right)^T &=&\Sigma_x^{-1} \left(X^{q(l+1)}-R^{q(l+1)-1} \right) \nonumber \\
									    &  + &  \left(\frac{\partial h}{\partial x} \mid _{x=X^{q(l+1)}}\right)^T
										\Sigma_{z}^{-1}  \left(h(X^{q(l+1)})-z^{l+1}\right). \nonumber \\
\end{eqnarray}
Here, we dropped the index $j$ for the particles for notational convenience.
We initialize the minimization using the result of a simplified implicit particle filter (see next subsection). Once the minimum is obtained, we substitute the random map (\ref{eq:RandomMap}) with $L_j=I$, where $I$ is the identity matrix, into (\ref{eq:AlgEq}) and solve the resulting scalar equation by Newton's method. The scalar derivative we need for the Newton steps is computed numerically. We initialize this iteration with $\lambda_j=0$. Finally, we compute the weights according to (\ref{eq:Weight}). If some weights are small, as indicated by a small effective sample size \cite{Arulampalam2002}, 
\begin{equation}
	\label{eq:EffectiveSampleSize}
	M_{Eff} = 1/\left(\displaystyle\sum\limits_{j=1}^M \left(w_j^{q(l+1)}\right)^2 \right)
\end{equation}
we resample using Algorithm 2 in \cite{Arulampalam2002}. The implicit filtering algorithm with gradient descent minimization and random maps is summarized in pseudo-code in Algorithm~\ref{algorithm:IPSFullNoise}.

\begin{algorithm}
\caption{Implicit Particle Filter with Random Maps and Gradient Descent Minimization}
\label{algorithm:IPSFullNoise}
\begin{algorithmic}
\STATE $ $
\STATE \COMMENT{Initialization, \ensuremath{t=0}}
\FOR{$j=1,\dots,M$}
	\STATE $\bullet$ sample $X_j^0 \sim p_o(X)$
\ENDFOR

\STATE $ $
\STATE $ $
\COMMENT{Assimilate observation \ensuremath{z^l}}
\FOR{$j=1,\dots,M$}
	   \STATE $\bullet$ Set up and minimize $F_j$ using gradient descent to compute $\phi_j$ and $\mu_j$
	   \STATE $\bullet$ Sample reference density $\xi_j\sim \mathcal{N}(0,I)$ 
	   \STATE $\bullet$ Compute $\rho_j=\xi_j^T\xi_j$ and $\eta_j = \xi_j/ \sqrt{\rho_j}$
	   \STATE $\bullet$ Solve (\ref{eq:AlgEq}) using the random map (\ref{eq:RandomMap}) with $L_j=I$
	   	   \STATE $\bullet$ Compute weight of the particle using (\ref{eq:Weight})
	      \STATE $\bullet$ Save particle $X_j$ and weight $w_j$
\ENDFOR
\STATE $ $
\STATE $\bullet$ Normalize the weights so that their sum equals 1
\STATE $\bullet$ Compute state estimate from $X_j$ weighted with $w_j$ (e.g. the mean)
\STATE $\bullet$ Resample if $M_{Eff}<c$
\STATE $\bullet$ Assimilate $z^{l+1}$
\end{algorithmic}
\end{algorithm}

This implicit filtering algorithm shares with weak constraint 4D-Var that a ``cost function'' (here $F_j$) is minimized by gradient descent. However, most 4D-Var implementations use tangent linear adjoint equations to compute the gradient. In the implicit filtering Algorithm \ref{algorithm:IPSFullNoise}, we do a fully nonlinear calculation of the gradient. Two further differences between 4D-Var and Algorithm \ref{algorithm:IPSFullNoise} are (\emph{i}) 4D-Var does not update the state sequentially, but the implicit particle filter does and, thus, reduces memory requirements; (\emph{ii}) 4D-Var computes the most likely state by minimizing the cost function, and this estimate can be biased; the implicit particle filter approximates the target density and, thus, can compute other statistics as state estimates, in particular the conditional expectation, which is, under wide conditions, the optimal state estimate (see e.g. \cite{ChorinHald}).

\subsection{A simplified implicit particle filtering algorithm with random maps and gradient descent minimization}\label{sec:SimplifiedFullNoise}
We wish to simplify the implicit particle filtering algorithm by reducing the dimension of the function $F_j$. The idea is to do an implicit sampling step only at times $t^{q(l+1)}$, i.e. when an observation becomes available. The state trajectory of each particle from time $t^{q(l)}$ (the last time an observation became available) to $t^{q(l+1)-1}$, is generated using the model equations (\ref{eq:DSDE}). This approach reduces the dimension of $F_j$ from $rm$ to $m$ (the state dimension). The simplification is thus very attractive if the number of steps between observations, $r$, is large. However, difficulties can also be expected for large $r$: the state trajectories up to time $t^{q(l+1)-1}$ are generated by the model alone and, thus, may not have a high probability with respect to the observations at time $t^{q(l+1)}$. The focussing effect of implicit sampling can be expected to be less emphasized and the number of particles required may grow as the gap between observations becomes larger. Whether or not the simplification we describe here can reduce the computational cost is problem dependent and we will illustrate advantages and disadvantages in the examples in Section \ref{sec:Geomagnetism}.

 Suppose we are given a collection of $M$ particles $X_j^{q(l)}$, $j=1,2,\dots,M$, whose empirical distribution approximates the target density at time $t^{q(l)}$ and the next observation, $z^{l+1}$, is available after $r$ steps at time $t^{q(l+1)}$. For each particle, we run the model for $r-1$ steps to obtain $X_j^{q(l)+1},\dots,X_j^{q(l+1)-1}$. We then define, for each particle, a function $F_j$ by
\begin{eqnarray}
	\label{eq:FFullNoiseDenseObs}
	F_j(X_j) &=&\frac{1}{2}  \left(X_j^{q(l+1)}-R_j^{q(l+1)-1}\right)^T \left(\Sigma_{x,j}^{q(l+1)-1}\right)^{-1}\left(X_j^{q(l+1)}-R_j^{q(l+1)-1}\right) \nonumber\\
			&+ &\frac{1}{2}  \left(h\left(X_j^{q(l+1)}\right)-z^{l+1}\right)^T \left(\Sigma_{z,j}^{q(l+1)}\right)^{-1}  \left(h\left(X_j^{q(l+1)}\right)-z^{l+1}\right),\nonumber \\
			&+&Z_j
\end{eqnarray}
whose gradient is given by (\ref{eq:LastGradFullNoise}). The algorithm then proceeds as Algorithm 1 in the previous section: we find the minimum of $F_j$ using gradient descent and solve (\ref{eq:AlgEq}) with the random map (\ref{eq:RandomMap}) with $L_j=I$. The weights are calculated by (\ref{eq:Weight}) with $r=1$ and the mean of $X_j$ weighted by $w_j$ is the state estimate at time $t^{q(l+1)}$.

This simplified implicit filter simplifies further if the observation function is linear, i.e. $h(x)=Hx$, where $H$ is a $k\times m$ matrix. One can show \cite{Morzfeld2011} that the minimim of $F_j$ is 
\begin{equation}
\label{eq:LinObsPhi}
	    \phi_j = \frac{1}{2}(z^{l+1}-HR_j^{q(l+1)-1})^T K_j^{-1}(z^{l+1}-HR_j^{q(l+1)-1}),
\end{equation}
with 
\begin{equation}
	  K_j =H\Sigma_{x,j}^{q(l+1)-1}H^T+\Sigma_{z,j}^{l+1}.
\end{equation}
A numerical approximation of the minimum is thus not required. The location of the minimum is 
\begin{equation}
\label{eq:LinObsMu}
	\mu_j = \Sigma_j \left(\left(\Sigma_{x,j}^{q(l+1)-1}\right)^{-1}R_j^{q(l+1)-1}+H^T(\Sigma_{z,j}^{q(l+1)})^{-1}z^{l+1}\right),
\end{equation}
with
\begin{equation}
\label{eq:LinObsSig}
	\Sigma_j^{-1} = \left(\Sigma_{x,j}^{q(l+1)-1}\right)^{-1}+H^T(\Sigma_{z,j}^{l+1})^{-1}H.
\end{equation}
Moreover, if $L_j$ is a Cholesky factor of $\Sigma_j$, then $X_j=\mu_j+L_j\xi_j$ solves (\ref{eq:AlgEq}) and the weights simplify to
\begin{equation}
	\label{eq:LinObsWeights}
	w_j^{n+1} \propto w_j^{n}\exp(-\phi_j) \left|\det L_j\right|.
\end{equation}
For the special case of a linear observation function and observations available at every model step ($r=1$), the simplified implicit filter is the full implicit filter and reduces to a version of optimal importance sampling \cite{Arulampalam2002,Bocquet2010,Morzfeld2011,Chorin2010}.
 
\section{Implicit particle filtering for equations with partial noise}
\label{sec:ImplicitSamplingPartialNoise}
We consider the case of a singular state covariance matrix $\Sigma_x$ in the context of implicit particle filtering. We start with an example taken from \cite{Jentzen2009}, to demonstrate how a singular state covariance appears naturally in the context of SPDE's driven by spatially smooth noise. The example serves as a motivation for more general developments in later sections. Another class of models with partial noise consists of dynamical equations supplemented by conservation laws. The dynamics are often uncertain and thus driven by noise processes, however there is typically zero uncertainty in the conservation laws (e.g. conservation of mass), so that the full model (dynamics and conservation laws) is subject to partial noise \cite{Kurapov2007}.

\subsection{Example of a model with partial noise: the semi-linear heat equation driven by spatially smooth noise}\label{sec:Heat}
We consider the stochastic semi-linear heat equation on the one-dimensional domain $x\in \left[0,1\right]$ over the time interval $t\in \left[0,1\right]$
\begin{equation}
	  \label{eq:SHE}
	\frac{\partial u}{\partial t} = \frac{\partial^2u}{\partial x^2}+\Gamma(u) + \frac{\partial W_t}{\partial t},
\end{equation}
where $\Gamma$ is a continuous function, and $W_t$ is a cylindrical Brownian motion (BM) \cite{Jentzen2009}. The derivative $\partial W_t/ \partial t$ in (\ref{eq:SHE}) is formal only (it does not exist in the usual sense). Equation (\ref{eq:SHE}) is supplemented by homogeneous Dirichlet boundary conditions and the initial value $u(x,0)=u_o(x)$. We expand the cylindrical BM $W_t$ in the eigenfunctions of the Laplace operator 
\begin{equation}	
	\label{eq:CBM}
	W_t = \displaystyle\sum\limits_{k=1}^{\infty} \sqrt{2q_k} \sin(k\pi x)\beta_t^k,
\end{equation}
where $\beta_t^k$ denote independent BM's and where the coefficients $q_k\geq 0$ must be chosen such that, for $\gamma \in (0,1)$,
\begin{equation}
	\label{eq:CBMInfSum}
	\displaystyle\sum\limits_{k=1}^{\infty} \lambda_k^{2\gamma-1}q_k<\infty,
\end{equation}
where $\lambda_k$ are the eigenvalues of the Laplace operator \cite{Jentzen2009}. If the coefficients $q_k$ decay fast enough, then, by (\ref{eq:CBM}) and basic properties of Fourier series, the noise is smooth in space and, in addition, the sum (\ref{eq:CBMInfSum}) remains finite as is required. For example one may be interested in problems where
\begin{equation}
\label{eq:Example_q_k}
	 q_k=\left\{ \begin{array}{l r} e^{-2k}, & \mbox{if }k\leq c, \\ 0,  & \mbox{if }k>c, \end{array}\right. 
\end{equation}
for some $c>0$.

The continuous equation must be discretized for computations and here we consider the Galerkin projection of the SPDE into an $m$-dimensional space spanned by the first $m$ eigenfunctions $e_k$ of the Laplace operator
\begin{equation}
	dU_t^m = (A_mU_t^m+\Gamma_m(U_t^m))dt+dW_t^m,
\end{equation}
where $U_t^m$, $\Gamma_m$ and $W_t^m$ are $m$-dimensional truncations of the solution, the function $\Gamma$ and the cylindrical BM $W_t$ respectively, and where $A_m$ is a discretization of the Laplace operator. Specifically, from (\ref{eq:CBM}) and (\ref{eq:Example_q_k}), we obtain:
\begin{equation}
		\label{eq:dWtExample}
	dW_t^m =  \displaystyle\sum\limits_{k=1}^{c} \sqrt{2}e^{-k} \sin(k\pi x)d\beta_t^k.
\end{equation}
After multiplying with the basis functions and integrating over the spatial domain, we are left with a set of $m$ stochastic ordinary differential equations
\begin{equation}
	\label{eq:Example}
	dx = f(x)dt+gdW
\end{equation}
where $x$ is an $m$-dimensional state vector, $f$ is a nonlinear vector function, $W$ is a BM. In particular, we calculate from (\ref{eq:dWtExample}):
\begin{equation}
	g = \frac{1}{\sqrt{2}}\mbox{diag}\left(\left(e^{-1},e^{-2},\dots,e^{-c},0,0,\dots,0\right)\right), \quad c<m,
\end{equation}
where $\mbox{diag}(a)$ is a diagonal matrix whose diagonal elements are the components of the vector $a$. Upon time discretization using, for example, a stochastic version of forward Euler with time step $\delta$ \cite{KloedenPlaten}, we arrive at (\ref{eq:DSDE}) with
\begin{equation}
	R(x) = x^n+\delta f(x^n), \quad
	  G(x) = \sqrt{\delta}g.
\end{equation}
It is now clear that the state covariance matrix $\Sigma_x = GG^T$ is singular for $c<m$. 

A singular state covariance causes no problems for running the discrete time model (\ref{eq:DSDE}) forward in time. However problems do arise if we want to know the pdf of the current state given the previous one. For example, the functions $F_j$ in the implicit particle filter algorithms (either those in Section \ref{sec:ImplicitSampling}, or those in \cite{Chorin2009,Chorin2010,Morzfeld2011}) are not defined for singular $\Sigma_x$. If $c\geq m$, then $\Sigma_x$ is ill-conditioned and causes a number of numerical issues in the implementation of these implicit particle filtering algorithms and, ultimately, the algorithms fail.

\subsection{Implicit particle filtering of models with partial noise, supplemented by densely available data}\label{sec:PartialNoiseDenseData}
We start deriving the implicit filter for models with partial noise by considering the special case in which observations are available at every model step ($r=1$). For simplicity, we assume that the noise is additive, i.e. $G(x^n,t^n) = G=\mbox{constant}$ and that $Q$ in (\ref{eq:Obs}) is a constant matrix. Under these assumptions, we can use a linear coordinate transformation to diagonalize the state covariance matrix and rewrite the model (\ref{eq:DSDE}) and the observations (\ref{eq:Obs}) as
\begin{eqnarray}
	\label{eq:forced}
	x^{n+1}&=&f(x^n,y^n,t^n) +  \Delta W^{n+1},\quad \Delta W^{n+1}\sim\mathcal{N}(0,\hat\Sigma_x) \\
	\label{eq:unforced}
	y^{n+1}&=&g(x^n,y^n,t^n), \\
	\label{eq:ObsPartial}
	z^{n+1}&=&h(x^n,y^n)+ QV^n,
\end{eqnarray}
where $x$ is a $p$-dimensional column vector, $p<m$ is the rank of the state covariance matrix (\ref{eq:StateCovariance}), and where $f$ is a $p$-dimensional vector function, $\hat\Sigma_x$ is a non-singular, diagonal $p\times p$ matrix, $y$ is a $(m-p)$-dimensional vector, and $g$ is a $(m-p)$-dimensional vector function. For ease of notation, we drop the hat above the ``new'' state covariance matrix $\hat\Sigma_x$ in (\ref{eq:forced}) and, for convenience, we refer to the set of variables $x$ and $y$ as the ``forced'' and ``unforced variables'' respectively. 

The key to filtering this system is observing that the unforced variables at time $t^{n+1}$, given the state at time $t^n$, are not random. To be sure, $y^n$ is random for any $n$ due to the nonlinear coupling $g(x,y)$, but the conditional pdf $p(y^{n+1}\mid x^n,y^n)$ is the delta-distribution. For a given (not random) initial state $x^0$, $y^0$, the target density is
 \begin{eqnarray}
 	\label{eq:PDFPartialNoiseDenseData}
	p(x^{0:n+1},y^{0:n+1}\mid z^{1:n+1})&\propto& p(x^{0:n},y^{0:n}\mid z^{1:n})  \nonumber \\
					&\times &p(z^{n+1}\mid x^{n+1},y^{n+1}) p(x^{n+1}\mid x^n,y^n). 
\end{eqnarray}

Suppose we are given a collection of $M$ particles, $X_j^{n},Y_j^{n}$, $j=1,2,\dots,M$, whose empirical distribution approximates the target density $p(x^{0:n},y^{0:n}\mid z^{1:n})$ at time $t^{n}$. The pdf for each particle at time $t^{n+1}$ is thus given by (\ref{eq:PDFPartialNoiseDenseData}) with the substitution of $X_j$ for $x$ and $Y_j$ for $y$. In agreement with the definition of $F_j$ in previous implementations of the implicit filter, we define $F_j$ here by
\begin{equation}
	\exp(-F_j(X_j^{n+1}))=p(z^{n+1}\mid X_j^{n+1},Y_j^{n+1}) p(X_j^{n+1}\mid X_j^n,Y_j^n).
\end{equation}
More specifically,
\begin{eqnarray}
	\label{eq:FPartialNoiseDenseObs}
	F_j(X_j^{n+1}) &=&\frac{1}{2}  \left(X_j^{n+1}-f_j^{n}\right)^T\Sigma_{x}^{-1}\left(X_j^{n+1}-f_j^{n}\right) \nonumber\\
			&+ &\frac{1}{2}  \left(h\left(X_j^{n+1},Y_j^{n+1}\right)-z^{n+1}\right)^T \Sigma_{z}^{-1}  \left(h\left(X_j^{n+1},Y_j\right)-z^{n+1}\right),\nonumber \\
			&+& Z_j
\end{eqnarray}
where $f_j^n$ is shorthand notation for $f(X_j^n,Y_j^n,t^n)$. With this $F_j$, we can use Algorithm \ref{algorithm:IPSFullNoise} to construct the implicit filter. For this algorithm we need the gradient of~$F_j$:
\begin{eqnarray}
	(\nabla F_j)^T &=& \Sigma_x^{-1} \left(X_j^{n+1}-f_j^{n} \right) \nonumber \\
									    &  + &  \left(\frac{\partial h}{\partial x} \mid _{x=X_j^{n+1}}\right)^T
										\Sigma_{z}^{-1}  \left(h(X_j^{n+1},Y_j^{n+1})-z^{n+1}\right). 
\end{eqnarray}
Note that $Y_j^{n+1}$ is fixed for each particle, if its previous state, $(X_j^{n},Y_j^{n})$, is known, so that the filter only updates $X_j^{n+1}$ when the observations $z^{n+1}$ become available. The unforced variables of the particles, $Y_j^{n+1}$, are moved forward in time using the model, as they should be, since there is no uncertainty in $y^{n+1}$ given $x^n,y^n$. The data are used in the state estimation of $y$ indirectly through the weights and through the nonlinear coupling between the forced and unforced variables of the model. If one observes only the unforced variables, i.e. $h(x,y)=h(y)$, then the data is not used directly when generating the forced variables, $X_j^{n+1}$, because the second term in (\ref{eq:FPartialNoiseDenseObs}) is merely a constant. In this case, the implicit filter becomes equivalent to a standard SIR filter, with weights $w_j^{n+1}=w_j^n \exp(-\phi_j)$.

This implementation of the implicit filter is numerically effective for filtering systems with partial noise, because the filter operates in a space of dimension $p$ (the rank of the state covariance matrix), which typically is less than the state dimension (see the example in Section \ref{sec:Geomagnetism}). The use of a gradient descent algorithm and random maps further makes the often costly computation of the Hessian of $F_j$ unnecessary. If $h$ is linear no iterative minimization is required.

If the state covariance matrix is ill-conditioned, a direct implementation of Algorithm \ref{algorithm:IPSFullNoise} is not possible. We propose to diagonalize the state covariance and set all eigenvalues below a certain threshold to zero so that a model of the form (\ref{eq:forced})-(\ref{eq:ObsPartial}) can be obtained. In our experience, such approximations are accurate and the filter of this section can be used.

\subsection{Implicit particle filtering for models with partial noise, supplemented by sparsely available data}\label{sec:PartialNoiseSparseData}
We extend the results of Section \ref{sec:PartialNoiseDenseData} to the more general case of observations that are sparse in time. Again, the key is to realize that $y^{n+1}$ is fixed given $x^n,y^n$. For simplicity, we assume additive noise and a constant $Q$ in (\ref{eq:Obs}), so that the target density becomes
 \begin{eqnarray}
 	\label{eq:PDFPartialNoiseSparseData}
	p(x^{0:q(l+1)},y^{0:q(l+1)}\mid z^{1:l+1})&\propto& p(x^{0:q(l)},y^{0:q(l)}\mid z^{1:l})  \nonumber \\
					&\times &p(z^{l+1}\mid x^{q(l+1)},y^{q(l+1)})\nonumber \\
					&\times& p(x^{q(l+1)}\mid x^{q(l+1)-1},y^{q(l+1)-1})\nonumber \\
						        &\times &p(x^{q(l+1)-1}\mid x^{q(l+1)-2},y^{q(l+1)-2})\nonumber \\
						        & &\vdots\nonumber \\
						        &\times &p(x^{q(l)+1}\mid x^{q(l)},y^{q(l)})\nonumber 
\end{eqnarray}
Given a collection of $M$ particles, $X_j^{n},Y_j^{n}$, $j=1,2,\dots,M$, whose empirical distribution approximates the target density $p(x^{0:q(l)},y^{0:q(l)}\mid z^{1:l})$ at time $t^{q(l)}$, we define, for each particle, the function $F_j$ by
\begin{eqnarray}
	\exp(-F_j(X_j)) &=& p(z^{l+1}\mid X_j^{q(l+1)},Y_j^{q(l+1)}) \nonumber \\
					&\times& p(X_j^{q(l+1)}\mid X_j^{q(l+1)-1},Y_j^{q(l+1)-1})\nonumber \\
					& & \vdots \nonumber \\ 
					&\times& p(X_j^{q(l)+1}\mid X_j^{q(l)},Y_j^{q(l)})
\end{eqnarray}
where $X_j$ is shorthand for $X_j^{q(l)+1,\dots,q(l+1)}$, so that
 \begin{eqnarray}
	   \label{eq:FPartialNoise}
	F_j(X_j) &=&\frac{1}{2} \left(X_j^{q(l)+1}-f_j^{q(l)}\right)^T \Sigma_{x}^{-1}\left(X_j^{q(l)+1}-f_j^{q(l)}\right) \nonumber\\
			 & +&\frac{1}{2} \left(X_j^{q(l)+2}-f_j^{q(l)+1}\right)^T \Sigma_{x}^{-1}\left(X_j^{q(l)+2}-f_j^{q(l)+1}\right) \nonumber\\
			 & &  \vdots\nonumber \\
			&+ &\frac{1}{2}  \left(X_j^{q(l+1)}-f_j^{q(l+1)-1}\right)^T \Sigma_{x}^{-1}\left(X_j^{q(l+1)}-f_j^{q(l+1)-1}\right) \nonumber\\
			&+ &\frac{1}{2}  \left(h\left(X_j^{q(l+1)},Y_j^{q(l+1)}\right)-z^{l+1}\right)^T \Sigma_{z}^{-1} \nonumber \\
			& & \times \left(h\left(X_j^{q(l+1)},Y_j^{q(l+1)}\right)-z^{l+1}\right)+Z_j.
\end{eqnarray}
At each model step, the unforced variables of each particle depend on the forced and unforced variables of the particle at the previous time step, so that $Y_j^{q(l+1)}$ is a function of $X_j^{q(l)},X_j^{q(l)+1},\dots,X_j^{q(l+1)-1}$ and $f_j^{q(l+1)}$ is a function of $X_j^{q(l)+1},X_j^{q(l)+2},$ $\dots,X_j^{q(l+1)}$. The function $F_j$ thus depends on the forced variables only. However, the appearances of the unforced variables in $F_j$ make it rather difficult to compute derivatives. The implicit filter with gradient descent minimization and random maps (see Algorithm \ref{algorithm:IPSFullNoise}) is thus a good filter for this problem, because it only requires computation of the first derivatives of $F_j$, while other implementations (see \cite{Chorin2010,Morzfeld2011}) require second derivatives as well. 

The gradient of $F_j$ is given by the $rp$-dimensional row vector
\begin{equation}
\label{eq:GradPartialNoise}
	\nabla F_j = \left(\frac{\partial F_j}{\partial X_j^{q(l)+1}},\frac{\partial F_j}{\partial X_j^{q(l)+2}},\dots, \frac{\partial F_j}{\partial X_j^{q(l+1)}}\right)
\end{equation}
with
\begin{eqnarray}
\label{eq:Grad1}
	\frac{\partial F_j}{\partial X_j^k}^T &=& \Sigma_x^{-1} \left( X_j^k -f_j^{k-1}\right) \nonumber \\
	&+& \left(\frac{\partial f}{\partial x} \mid_{k}\right)^T \Sigma_x^{-1} \left( X_j^{k+1} -f_j^{k}\right) \nonumber \\
	   &+&\left(\frac{\partial f}{\partial y}\mid_{k+1}\frac{\partial y^{k+1}}{\partial X_j^k}\right)^T \Sigma_x^{-1} \left( X_j^{k+2} -f_j^{k+1}\right) \nonumber \\
	&+&\left(\frac{\partial f}{\partial y}\mid_{k+2}\frac{\partial y^{k+2}}{\partial X_j^k}\right)^T \Sigma_x^{-1} \left( X_j^{k+3} -f_j^{k+2}\right) \nonumber \\
	   & &\vdots\nonumber \\
	&+&\left(\frac{\partial f}{\partial y}\mid_{q(l)-1}\frac{\partial y^{q(l)-1}}{\partial X_j^k}\right)^T \Sigma_x^{-1} \left( X_j^{q(l+1)} -f_j^{q(l)-1}\right) \nonumber \\
	&+& \left(\frac{\partial h}{\partial y}\mid_{k}\frac{\partial y^{q(l)}}{\partial X_j^k}\mid_{k-1}\right)^T \Sigma_z^{-1} \left( h\left(X_j^{q(l+1)},Y_j^{q(l+1)}\right) -z^{l+1}\right) 
\end{eqnarray}
for $k=q(l)+1,\dots,q(l+1)-1$ and where $(\cdot)\mid_k$ denotes ``evaluate at time $t^k$.'' The derivatives $\partial y^{i}/\partial X_j^k$, $i=k+1,\dots,q(l)$, can be computed recursively while constructing the sum, starting with
\begin{equation}
\label{eq:Grad2}
	\frac{\partial y^{k+1}}{\partial X_j^k}= \frac{\partial }{\partial X_j^k}\left(g(X_j^k,Y_j^k)\right)=\frac{\partial g}{\partial x}\mid_k,
\end{equation}
and then using
\begin{equation}
	\label{eq:Grad3}
	\frac{\partial y^{k+i}}{\partial X_j^k}=\frac{\partial g}{\partial x}\mid_{i-1} \frac{\partial y^{i-1}}{\partial X_j^k}\mid_{i-1}\mbox{ },\quad i=k+2,\dots,q(l).
\end{equation}

The minimization of $F_j$ for each particle is initialized with a free model run for $r$ steps, with initial conditions given by the final position of the $j^{th}$ particle at the previous assimilation step. With this initial guess we compute the gradient using (\ref{eq:GradPartialNoise})-(\ref{eq:Grad3}) and, after a line search and one step of gradient descent, obtain a new set of forced variables. We use this result to update the unforced variables by the model, and proceed to the next iteration. Once the minimum $\phi_j$ and its location $\mu_j$ are found, we use the random map (\ref{eq:RandomMap}) with $L_j=I$ to compute $X_j^{q(l)+1},\dots,X_j^{q(l+1)}$ for this particle and then use these forced variables to compute $Y_j^{q(l)+1,\dots,q(l+1)}$. We do this for all particles, and compute the weights from (\ref{eq:Weight}) with $m=p$, then normalize the weights so that their sum equals one and thereby obtain an approximation of the target density. We resample if the effective sample size $M_{Eff}$ is below a threshold and move on to assimilate the next observation. The implicit filtering algorithm is summarized with pseudo code in Algorithm~\ref{algorithm:IPSPartialNoise}.
\begin{algorithm}
\caption{Implicit Particle Filter with Random Maps and Gradient Descent Minimization for Models with Partial Noise}
\label{algorithm:IPSPartialNoise}
\begin{algorithmic}
\STATE $ $
\STATE \COMMENT{Initialization, \ensuremath{t=0}}
\FOR{$j=1,\dots,M$}
	\STATE $\bullet$ sample $X_j^0 \sim p_o(X)$
\ENDFOR

\STATE $ $
\STATE $ $
\COMMENT{Assimilate observation \ensuremath{z^l}}

\FOR{$j=1,\dots,M$}
	\STATE $\bullet$ Set up and minimize $F_j$ using gradient descent:
	\STATE Initialize minimization with a free model run
	  \WHILE{Convergence criteria not satisfied}
 		\STATE Compute gradient by (\ref{eq:GradPartialNoise})
			\STATE Do a line search
			\STATE Compute next iterate by gradient descent step 
			\STATE Use results to update unforced variables using the model 			\STATE Check if convergence criteria are satisfied
		\ENDWHILE
	   \STATE $\bullet$ Sample reference density $\xi_j\sim \mathcal{N}(0,I)$ 
	   \STATE $\bullet$ Compute $\rho_j=\xi_j^T\xi_j$ and $\eta_j = \xi_j/\sqrt{\rho_j}$
	   \STATE $\bullet$ Solve (\ref{eq:AlgEq}) using random map (\ref{eq:RandomMap}) with $L_j=I$ to compute $X_j$
	   \STATE $\bullet$ Use this $X_j$ and the model to compute corresponding $Y_j$
	   \STATE $\bullet$ Compute weight of the particle using (\ref{eq:Weight})
	      \STATE $\bullet$ Save particle $(X_j,Y_j)$ and weight $w_j$
\ENDFOR
\STATE $ $
\STATE $\bullet$ Normalize the weights so that their sum equals 1
\STATE $\bullet$ Compute state estimate from $X_j$ weighted with $w_j$ (e.g. the mean)
\STATE $\bullet$ Resample if $M_{Eff}<c$
\STATE $\bullet$ Assimilate $z^{l+1}$
\end{algorithmic}
\end{algorithm}

Note that all state variables are computed by using both the data and the model, regardless of which set of variables (the forced or unforced ones) is observed. The reason is that sparse observations induce a nonlinear coupling, through $f$ and $g$ in (\ref{eq:forced})-(\ref{eq:ObsPartial}), between the unforced and forced variables at the various model steps. It should also be noted that the function $F_j$ is a function of $rp$ variables (rather than $rm$), because the filter operates in the subspace of the forced variables. If the minimization is computationally too expensive, because $p$ or $r$ is extremely large, then one can easily adapt the ``simplified'' implicit particle filter of Section \ref{sec:SimplifiedFullNoise} to the situation of partial noise using the methods we have described above. If $h$ is nonlinear, this simplified filter requires a minimization of a $p$-dimensional function for each particle. If $h$ is linear, no numerical minimization is required.

\subsection{Discussion}
We wish to point out similarities and differences between the implicit filter and three other data assimilation methods. In particular, we discuss how data are used in the generation of the state estimates. It is clear that the implicit filter uses the available data as well as the model to generate the state trajectories for each particle, i.e. it makes use of the nonlinear coupling between forced and unforced variables.

SIR and EnKF make less direct use of the data. In SIR, the particle trajectories are generated using the model alone and only later weighted by the observations. Data thus propagate to the SIR state estimates indirectly through the weights. In EnKF, the state trajectories are generated using the model and only the states at times $t^{q(l)}$ (when data are available) are updated by data. Thus, EnKF uses the data only to update its state estimates at times for which data are actually available. 

A weak constraint 4D-Var method is perhaps closest in spirit to the implicit filter. In weak constraint 4D-Var, a cost function similar to $F_j$ is minimized (typically by gradient descent) to find the state trajectory with maximum probability given data and model. If one picks the time window for a 4D-Var assimilation from one observation $z^l$ to the next $z^{l+1}$, then the use of the data is similar to the use of the data in an implicit filter, because, in both algorithms, the model as well as data are used to generate the state trajectories. In fact, one can view the implicit particle filter as a randomized and sequential version of weak constraint 4D-Var (or, one may interpret weak constraint 4D-Var as an implicit smoother with a single particle). These issues will be taken up in more detail in a subsequent paper \cite{Atkins2011}. 

Finally, we would like to discuss the implicit filtering algorithm for the special case of a perfect model, i.e. 
\begin{eqnarray}
	\label{eq:perfectModel1}
	y^{n+1}&=&g(y^n,t^n), \\
	\label{eq:perfectModel2}
	z^{l}&=&h(y^{q(l)})+ Q^{l}V^{l}.
\end{eqnarray}
Following the steps above and, assuming we are given a collection of $M$ particles, $Y_j^{n}$, $j=1,2,\dots,M$, whose empirical distribution approximates the target density $p(y^{0:q(l)}\mid z^{1:l})$ at time $t^{q(l)}$, we define, for each particle, the function $F_j$ by
\begin{equation}
	\exp(-F_j) = p(z^{l+1}\mid Y_j^{q(l+1)}),
\end{equation}
so that
 \begin{eqnarray}
	F_j &=&\frac{1}{2}  \left(h\left(Y_j^{q(l+1)}\right)-z^{l+1}\right)^T \left(\Sigma_{z,j}^{q(l+1)}\right)^{-1}  \times \left(h\left(Y_j^{q(l+1)}\right)-z^{l+1}\right) \nonumber \\
			&+&Z_j
\end{eqnarray}
Since $Y_j^{q(l)}$ is fixed for each particle, $F_j$ is merely a constant that is used to weigh the particle trajectory by the weight
\begin{equation}
	w_j^{l+1}=w_j^l \exp(-F_j).
\end{equation}
The data are used indirectly here, because the initial condition determines the full state trajectory. However, this initial condition is fixed for each particle. For a perfect model, strong constraint 4D-Var makes more efficient use of the available data by using it to find an ``optimal initial condition,'' compatible with the data.

\section{Application to Geomagnetism}
\label{sec:Geomagnetism}
\begin{figure}[htbp]
\begin{center}
{\includegraphics[width=1\textwidth]{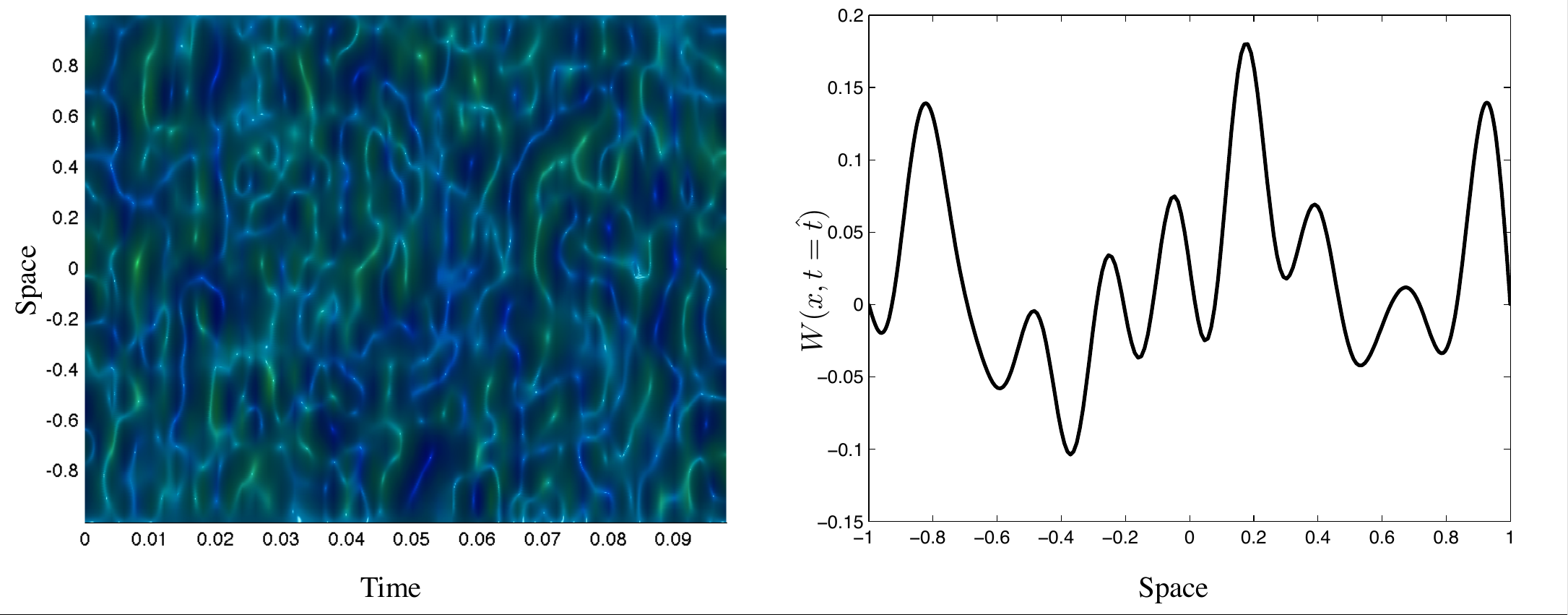}}
\caption{The noise process $W(x,t)$. Left: Noise plotted as a function of $x$ and $t$. Right: Snapshot of $W(x,t=\hat t)$.} 
\label{fig:noise} 
\end{center}
\end{figure}

Data assimilation has been recently applied to geomagnetic applications and there is a need to find out which data assimilation technique is most suitable \cite{Fournier2010}. Thus far, a strong constraint 4D-Var approach \cite{Fournier2007} and a Kalman filter approach \cite{Sun2007} have been considered. Here, we apply the implicit particle filter to a test problem very similar to the one first introduced by Fournier and his colleagues in \cite{Fournier2007}. The model is given by two SPDE's
\begin{eqnarray}
	\label{eq:u}
	\partial_tu+u\partial_xu &=& b\partial_xb+\nu \partial_x^2u + g_u \partial_t W(x,t), \\
	\label{eq:b}
	\partial_tb+u\partial_xb&=&b\partial_xu+\partial_x^2b +g_b\partial_t W(x,t),
\end{eqnarray}
where, $g_u,g_b$ are scalars, and where $W$ is a stochastic process (the derivative here is formal and may not exist in the usual sense). We study the above equations with $\nu =10^{-3}$ as in \cite{Fournier2007}, and with $g_u=0.01$, $g_b=1$, so that the uncertainty in the unobserved quantity is much larger than the uncertainty in the unobserved quantity. We consider the above equations on the strip $0\leq t \leq T$, $-1\leq x\leq1$ and with boundary conditions
\begin{eqnarray}
	\label{eq:InitialU}
	u(x,t)=0, \mbox{   if } x=\pm 1, & & u(x,0)=\sin(\pi x)+2/5\sin(5\pi x), \\
	\label{eq:InitialB}
	b(x,t)=\pm1, \mbox{   if } x=\pm 1, & & b(x,0)=\cos(\pi x)+2\sin(\pi (x+1)/4).
\end{eqnarray}
Physically, $u$ represents the velocity field and $b$ represents the secular variation of the magnetic field. The model is essentially the model proposed in \cite{Fournier2007}, but with additive noise
\begin{equation}
	\label{eq:NoiseProcess}
	W(x,t)=\sum_{k=0}^{\infty} \alpha_k \sin(k\pi x) w_k^1(t) + \beta_k \cos(k \pi/2x) w_k^2(t).
\end{equation}
where $w^1_k,w^2_k$ are independent BMs and where
\begin{equation}
	\label{eq:ChoiceCoeff}
	\alpha_k = \beta_k = \left\{ \begin{array}{l r} 1, & \mbox{if } k\leq 10, \\ 0,  & \mbox{if }k>10. \end{array}\right.
\end{equation}
Here, we are content with this simple noise model that represents a small uncertainty at the boundaries of both fields and is spatially smooth. However, it is straightforward to incorporate more information about the spatial distribution of the uncertainty by picking suitable coefficients $\alpha_k$, $\beta_k$. An illustration of the noise process is given in Figure \ref{fig:noise}. 

\subsection{Discretization of the dynamical equations}
We follow \cite{Fournier2007} in the discretization of the dynamical equations, however we decided to present some details of the discretization to explain how the noise process $W$ comes into play.

For both fields, we use Legendre spectral elements of order $N$ (see e.g. \cite{Canuto2006,Deville2006}), so that
\begin{eqnarray}
	u(x,t)&=&\displaystyle\sum\limits_{j=0}^N\hat{u}_j(t)\psi_j(x)=\displaystyle\sum\limits_{j=1}^{N-1}\hat{u}_j(t)\psi_j(x),\nonumber \\
	   b(x,t)&=&\displaystyle\sum\limits_{j=0}^N\hat{b}_j(t)\psi_j(x)=-\psi_1(x)+\psi_N(x)+ \displaystyle\sum\limits_{j=1}^{N-1}\hat{b}_j(t)\psi_j(x), \nonumber \\
	   W(x,t)&=&\displaystyle\sum\limits_{j=0}^{N} \hat{W}_{j}(t)\psi_j(x)=\displaystyle\sum\limits_{j=1}^{N-1} \hat{W}_{j}(t)\psi_j(x)  \nonumber
\end{eqnarray}
where $\psi_j$ are the characteristic Lagrange polynomials of order $N$, centered at the $j$th Gauss-Lobatto-Legendre (GLL) node $\xi_j$. We substitute the expansions into the weak form of (\ref{eq:u}) and (\ref{eq:b}) (no integration by parts) and evaluate the integrals by Gauss-Lobatto-Legendre quadrature
\begin{equation}
	\int_{-1}^{1}p(x)dx \sim \displaystyle \sum_{j=0}^Np(\xi_j)w_j, \nonumber
\end{equation}
where $w_j$ are the corresponding weights. Making use of the orthogonality of the basis functions, $\psi_j(\xi_k)=\delta_{j,k}$, we obtain the set of SDE's
\begin{eqnarray}
	M\partial_t\hat u&=&M\left(\hat b\circ D\hat b-\hat u\circ D\hat u+\nu D^2\hat u+ \Psi_x^B \hat b +g_u\partial_t \hat W\right), \nonumber\\
	M\partial_t\hat b&=&M\left(\hat b\circ D\hat u-\hat u\circ D\hat b+ D^2\hat b- \Psi_x^B \hat u+\Psi_{xx}^B+g_b\partial_t \hat W  \right), \nonumber
\end{eqnarray}
where $\circ$ denotes the Hadamard product ($(\hat u\circ \hat b)_k=\hat u_k \hat b_k$), $\hat u, \hat b, \hat W$ are $m=(N-2)$-dimensional column vectors whose components are the coefficients in the series expansions of $u,b,W_u$ and $W_b$ respectively and where $\Psi_x^B=\mbox{diag}\left((\partial_x\psi_j(\xi_1),\dots,\partial_x\psi_j(\xi_{N-1}))\right)$ and $\Psi_{xx}^B=(\partial_{xx}\psi_2(\xi_1),\dots,\partial_{xx}\psi_{N-1}(\xi_{N-1}))^T$ is a diagonal $m \times m$ matrix and an $m$-dimensional column vector respectively, which make sure that our approximation satisfies the boundary conditions. In the above equations, the $m\times m$ matrices $M$, $D$ and $D^2$ are given by
\begin{eqnarray}
	M=\mbox{diag}\left((w_1,\dots,w_{N-1})\right),\quad D_{j,k}=\partial_x\psi_j(\xi_k),\quad D^2_{j,k}=\partial_{xx}\psi_j(\xi_k).\nonumber 
\end{eqnarray}
We apply a first-order implicit-explicit method with time step $\delta$ for time discretization and obtain the discrete-time and discrete-space equations
\begin{eqnarray}
	(M-\delta \nu MD^2)u^{n+1}&=&M\left(u^n+\delta\left(b^n\circ Db^n-u^n\circ Du^n+ \Psi_x^B b^n \right) \right)+\Delta W_u^n, \nonumber\\
	(M-\delta MD^2)b^{n+1}&=&M\left(b^n +\delta\left(b^n\circ Du^n-u^n\circ Db^n- \Psi_x^B u^n+\Psi_{xx}^B  \right)\right)+ \Delta W_b^n, \nonumber
\end{eqnarray}
where
\begin{equation}
		\Delta W_u \sim \mathcal{N}(0,\Sigma_u),\quad \Delta W_b \sim \mathcal{N}(0,\Sigma_b) 
\end{equation}
and
\begin{eqnarray}
\label{eq:SigmaU}
\Sigma_u &=&g_u^2 \delta M\left(
F_sCC^TF_s^T+F_cCC^TF_c^T
\right)M^T,  \\
\label{eq:SigmaB}
\Sigma_b &=&g_b^2\delta M\left(
F_sCC^TF_s^T+F_cCC^TF_c^T
\right)M^T,  \\
C& =& \mbox{diag}((\alpha_1,\dots,\alpha_n)),  \\
F_s &=& (\sin(\pi),\sin(2\pi),\dots,\sin(m\pi))(\xi_1, \xi_2,\dots,\xi_m)^T,  \\
F_c &=& (\cos(\pi/2),\cos(3\pi/2),\dots,\cos(m\pi/2))(\xi_1, \xi_2,\dots,\xi_m)^T. 
\end{eqnarray}

For our choice of $\alpha_k,\beta_k$ in (\ref{eq:ChoiceCoeff}), the state covariance matrices $\Sigma_u$ and $\Sigma_b$ are singular if $m>10$. To diagonalize the state covariances we solve the symmetric eigenvalue problems \cite{Parlett}
\begin{eqnarray}
(M-\delta \nu MD^2)v_u  = \Sigma_u v_u \lambda^u,\nonumber \\ 
(M-\delta  MD^2)v_b = \Sigma_b v_b\lambda^b,\nonumber 
\end{eqnarray}
and define the linear coordinate transformations
\begin{equation}
	u=V_u(x_u,y_u)^T,\quad b=V_b(x_b,y_b)^T,
\end{equation}
where the columns of the $m\times m$-matrices $V_u$ and $V_b$ are the eigenvectors of $v_u$, $v_b$ respectively. The discretization using Legendre spectral elements works in our favor here, because the matrices $M$ and $D^2$ are symmetric so that we can diagonalize the left hand side simultaneously with the state covariance matrix to obtain
\begin{eqnarray}
x_u^{n+1} &=& f_u(x_u^n,y_u^n,x_b^n,y_b^n)+\Delta \hat W_u^n, \nonumber \\
y_u^{n+1} &=& g_u(x_u^n,y_u^n,x_b^n,y_b^n), \nonumber \\
x_b^{n+1} &=& f_b(x_u^n,y_u^n,x_b^n,y_b^n)+\Delta \hat W_b^n, \nonumber \\
y_b^{n+1} &=& g_b(x_u^n,y_u^n,x_b^n,y_b^n), \nonumber 
\end{eqnarray}
where $f_u,f_b$ are 10-dimensional vector functions, $g_u,g_b$ are $(m-10)$-dimensional vector functions and where 
\begin{eqnarray}
\hat W_u^n&\sim& \mathcal{N}\left(0,\mbox{diag}\left(\left(\lambda^u_{1},\lambda^u_{2},\dots,\lambda^u_{10}\right)\right)\right),\nonumber \\
\hat W_b^n&\sim& \mathcal{N}\left(0,\mbox{diag}\left(\left(\lambda^b_{1},\lambda^b_{2},\dots,\lambda^b_{10}\right)\right)\right).\nonumber
\end{eqnarray}

\begin{figure}[htbp]
\begin{center}
{\includegraphics[width=1\textwidth]{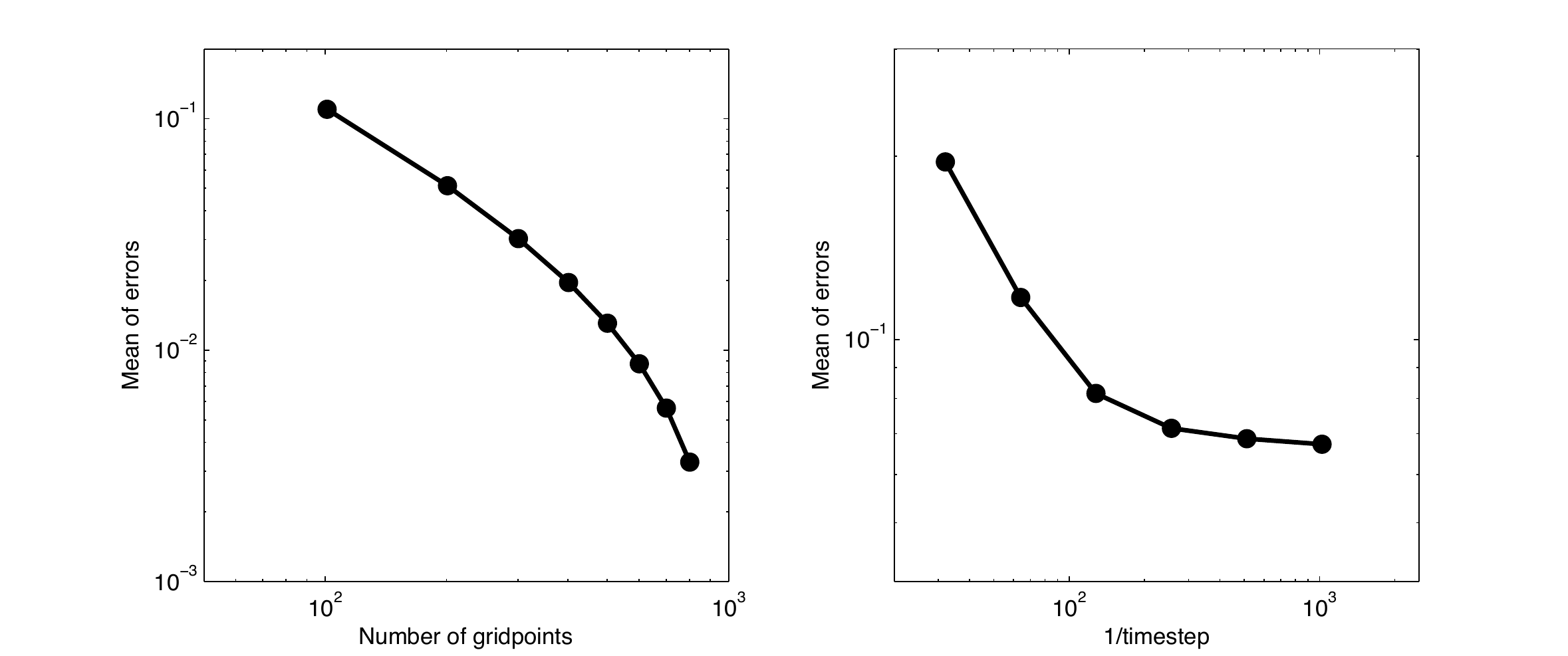}}
\caption{Convergence of discretization scheme for geomagnetic equations. Left: Convergence in the number of spatial grid-points. Right: Convergence in the time step.} 
\label{fig:Convergence} 
\end{center}
\end{figure}

We test the convergence of our approximation as follows. To assess the convergence in the number of grid-points in space, we define a reference solution using $N=2000$ grid-points and a time step of $\delta = 0.002$. We compute another approximation of the solution, using the same (discrete) BM as in the reference solution, but with another number of grid-points, say $N=500$. We compute the error at $t=T=0.2$, $e_x=\left|\right|(u_{500}(x,T)^T,b_{500}(x,T)^T)-(u_{Ref}(x,T)^T,b_{Ref}(x,T)^T)\left|\right|$, where $\left|\right|\cdot \left|\right|$ denotes the Euclidean norm, and store it. We repeat this procedure 500 times and compute the mean of the error norms. The results are shown in the left panel of Figure \ref{fig:Convergence}. 
We observe a super algebraic convergence as expected from a spectral method. 

Similarly, we check the convergence of the approximation in the time step by computing a reference solution with $N_{Ref}=1000$ and $\delta_{Ref} = 2^{-12}$. Using the same BM as in the reference solution, we compute an approximation with time step $\delta$ and compute the error at $t=T=0.2$, $e_t=\left|\right|(u_{\delta}(x,T)^T,b_{\delta}(x,T)^T)-(u_{Ref}(x,T)^T,b_{Ref}(x,T)^T)\left|\right|$, and store it. We repeat this procedure 500 times and then compute the mean of these error norms. The results are shown in the right panel of Figure \ref{fig:Convergence}. We observe a first order decay in the error as is expected. The scheme has converged for time steps smaller than $\delta = 0.002$, so that a higher resolution in time does not improve the accuracy of the approximation. Moreover, the Courant-Friederichs-Lewy condition limits the time step for a given number of nodes. The limit here is quite strict because the Legendre elements accumulate grid-points close to the boundaries so that the smallest spacing between grid-points is very small, even for a moderate number of nodes.

Here we are satisfied with an approximation with $\delta = 0.002$ and $N=300$ grid-points in space as in \cite{Fournier2007}. The relatively small number of spatial grid-points is sufficient because the noise is very smooth in space and because the Legendre spectral elements accumulate nodes close to the boundaries and, thus, represent the steep boundary layer, characteristic of (\ref{eq:u})-(\ref{eq:b}), well even if $N$ is small (see also \cite{Fournier2007}). 
 
\subsection{Filtering results}
\begin{figure}[htbp]
\begin{center}
{\includegraphics[width=1\textwidth]{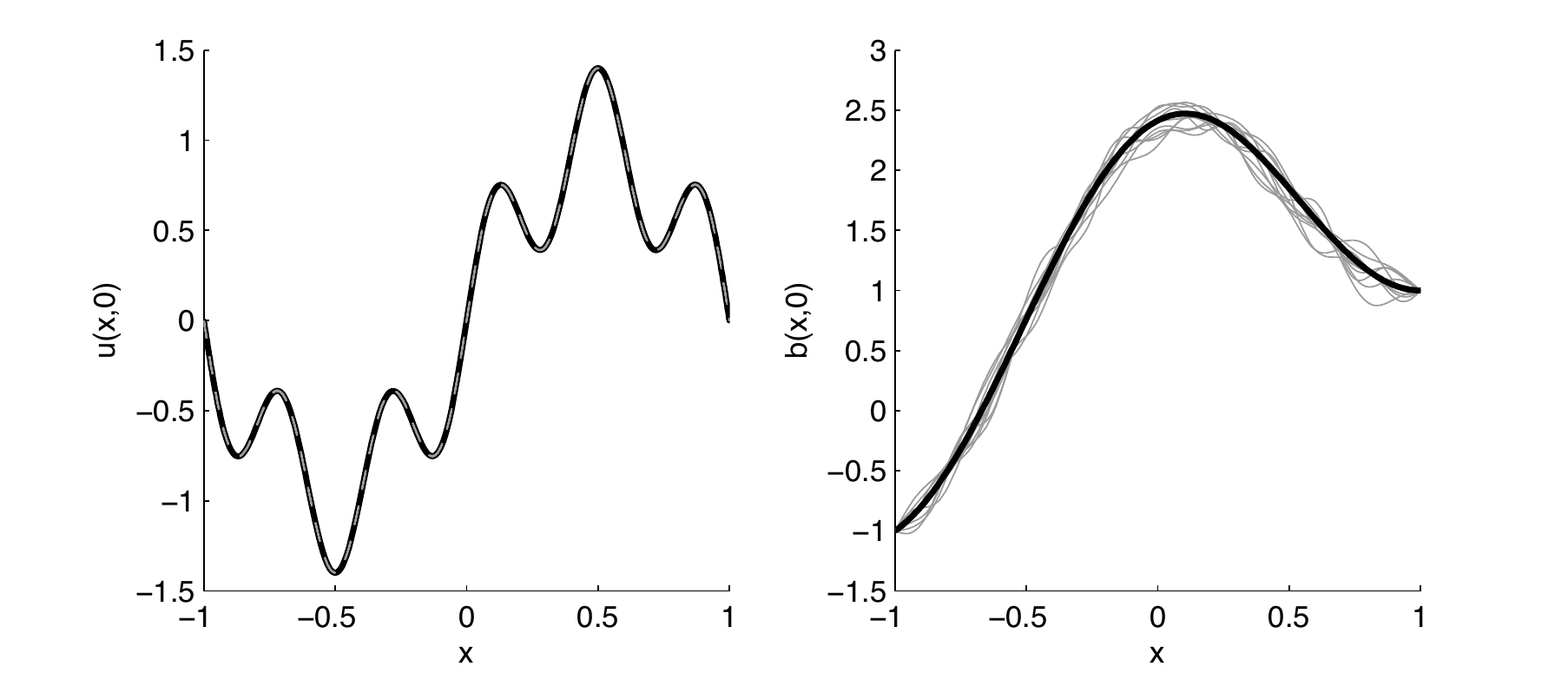}}
\caption{Uncertainty in the initial state. Left: $u(x,0)$ (unobserved). Right: $b(x,0)$ (observed). Black: mean. Grey: 10 realizations of the initial state.} 
\label{fig:Initial} 
\end{center}
\end{figure}

We apply the implicit particle filter with gradient descent minimization and random maps (see Algorithm \ref{algorithm:IPSPartialNoise} in Section \ref{sec:ImplicitSamplingPartialNoise}), the simplified implicit particle filter (see Section \ref{sec:SimplifiedFullNoise}) adapted to models with partial noise, a standard EnKF (without localization or inflation), as well as a standard SIR filter to the test problem (\ref{eq:u})-(\ref{eq:b}). The numerical model is given by the discretization described in the previous section with a random initial state. The distribution of the initial state is Gaussian with mean $u(x,0),b(x,0)$ as in (\ref{eq:InitialU})-(\ref{eq:InitialB}) and with a covariance $\Sigma_u,\Sigma_b$ given by (\ref{eq:SigmaU})-(\ref{eq:SigmaB}). In Figure \ref{fig:Initial}, we illustrate the uncertainty in the initial state and
plot 10 realizations of the initial state (grey lines) along with its mean (black lines). We observe that the uncertainty in $u_0$ is small compared to the uncertainty in $b_0$. 

The data are the values of the magnetic field $b$, measured at $k$ equally spaced locations in $\left[ 0,1\right]$ and corrupted by noise: 
\begin{equation}
	z^{l} = Hb^{q(l)}+sV^{l},
\end{equation}
where $s=0.001$ and where $H$ is a $k\times m$-matrix that maps the numerical approximation $b$ (defined at the GLL nodes) to the locations where data is collected. We consider data that are dense in time ($r=1$) as well as sparse in time ($r>1$). The data are sparse in space and we consider two cases: (\emph{i}) we collect the magnetic field at 200 equally spaced locations; and (\emph{ii}) we collect $b$ at 20 equally spaced locations. The variables $u$ are unobserved and it is of interest to study how the various data assimilation techniques make use of the information in $b$ to update the unobserved variables $u$ \cite{Fournier2007,Fournier2010}.

\begin{figure}[htbp]
\begin{center}
{\includegraphics[width=1\textwidth]{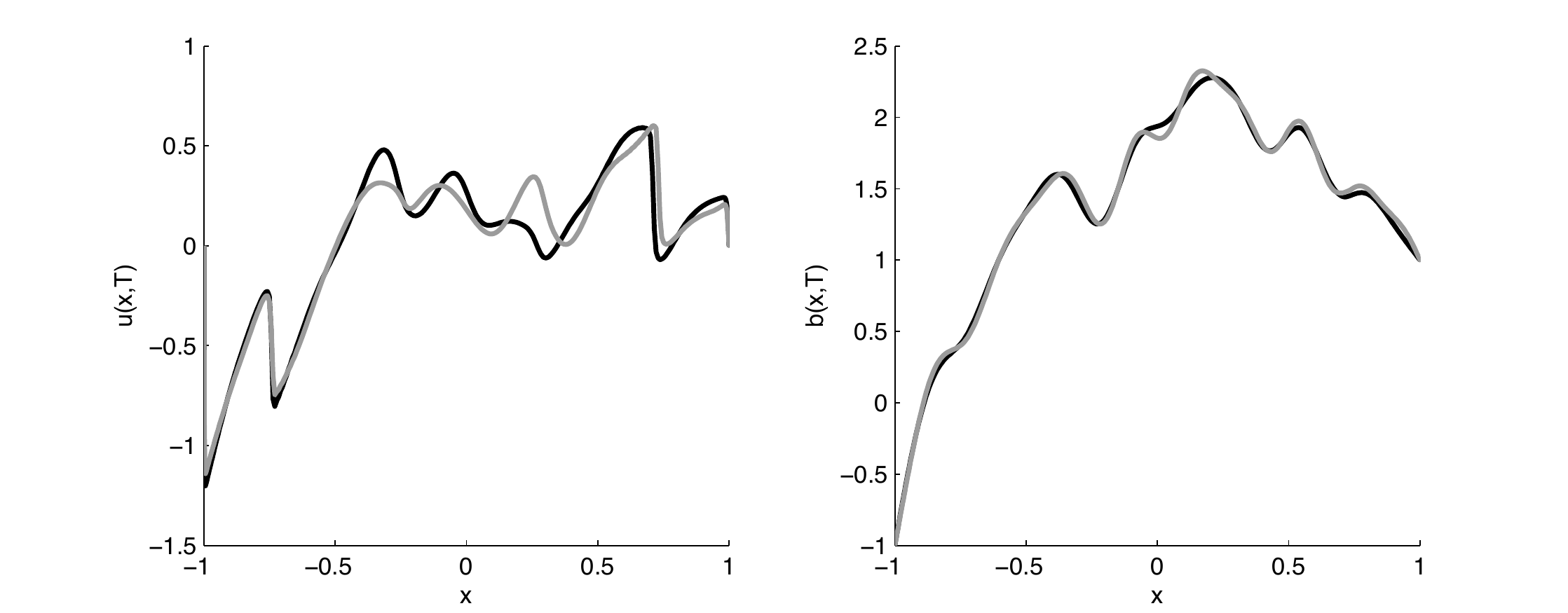}}
\caption{Outcome of a twin experiment. Black: true state $u(x,T)$ (left) and $b(x,T)$ (right). Grey: reconstruction by implicit particle filter with 4 particles.} 
\label{fig:Twin} 
\end{center}
\end{figure}

To assess the performance of the filters, we ran 100 twin experiments. A twin experiment amounts to: (\emph{i}) drawing a sample from the initial state and running the model forward in time until $T=0.2$ (one fifth of a magnetic diffusion time \cite{Fournier2007}) (\emph{ii}) collecting the data from this free model run; and (\emph{iii}) using the data as the input to a filter and reconstructing the state trajectory. Figure \ref{fig:Twin} shows the result of one twin experiment for $r=4$.

For each twin experiment, we calculate and store the error at $T=0.2$ in the velocity, $e_u = \left|\right|u(x,T)-u_{Filter}(x,T)\left|\right|/\left|\right|u(x,T)\left|\right|$, and in the magnetic field, $e_b = \left|\right|b(x,T)-b_{Filter}(x,T)\left|\right|/\left|\right|b(x,T)\left|\right|$. After running the 100 twin experiments, we calculate the mean of the error norms (not the mean error) and the variance of the error norms (not the variance of the error). All filters we tested were ``untuned,'' i.e. we have not adjusted or inserted any free parameters to boost the performance of the filters. 

\begin{figure}[htbp]
\begin{center}
{\includegraphics[width=0.6\textwidth]{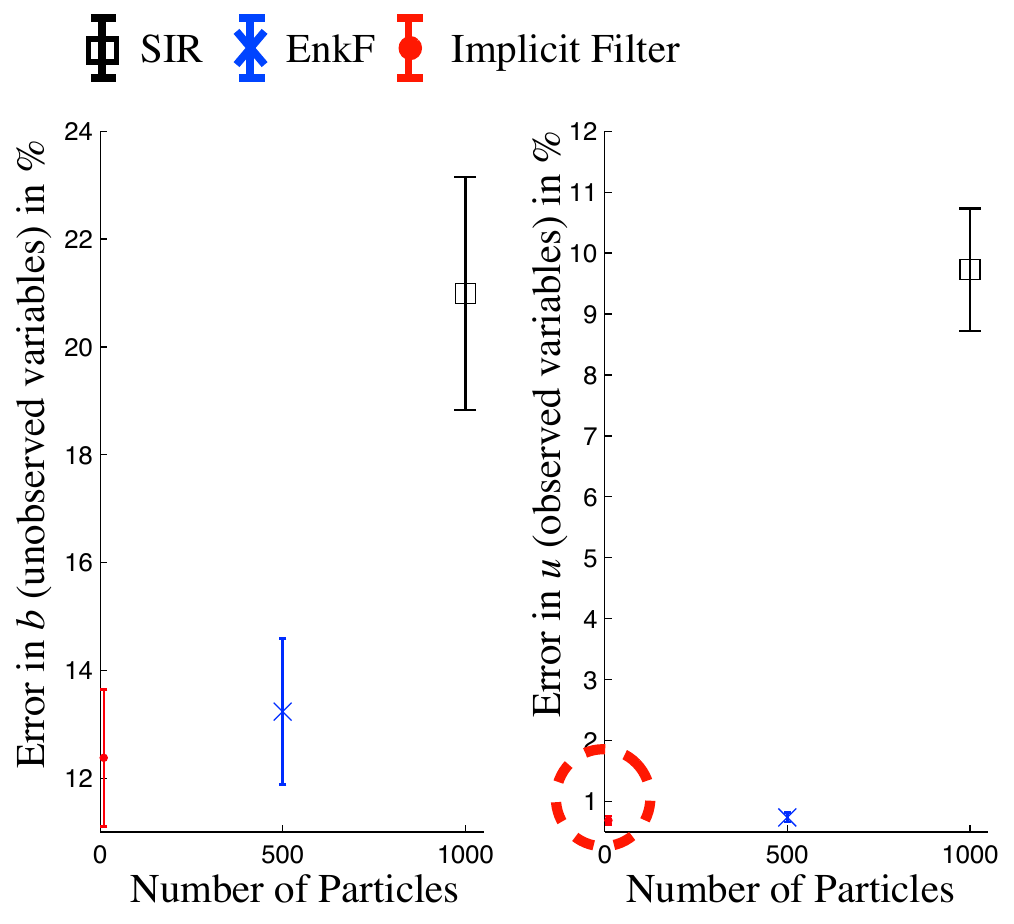}}
\caption{Filtering results for data collected at a high spatial resolution (200 measurement locations). The errors at $T=0.2$ of the implicit particle filter (red),  EnKF (blue) and SIR filter (black) are plotted as a function of the number of particles. The error bars represent the mean of the errors and mean of the standard deviations of the errors.} 
\label{fig:ResultsDense1} 
\end{center}
\end{figure}

Figure \ref{fig:ResultsDense1} shows the results for the implicit particle filter, the EnKF as well as the SIR filter for $200$ measurement locations and for $r=10$.
It is evident from this figure that the implicit particle filter requires only very few particles to yield accurate state estimates with less than 1\% error in the observed variables. The SIR filter with 1000 particles gives significantly larger errors (about 10\% in the observed variables) and much larger variances in the errors. The EnKF requires about 500 particles to come close to the accuracy of the implicit filter with only 10 particles.

In the experiments, we observed that the minimization in implicit particle filtering typically converged after 4-10 steps (depending on $r$, the gap in time between observations). The convergence criterion was to stop the iteration when the change in $F_j$ was less than $10\%$. A more accurate minimization did not improve the results significantly, so that we were satisfied with a relatively crude estimate of the minimum in exchange for a speed-up of the algorithm. We found $\lambda$ by solving (\ref{eq:AlgEq}) with Newton's method using $\lambda^0 = 0$ as initial guess and observed that it converged after about eight steps. The convergence criterion was to stop the iteration if $\left| F(\lambda)-\phi-\rho \right| \leq 10^{-3}$, because the accurate solution of this scalar equation is numerically inexpensive. We resampled using Algorithm 2 in \cite{Arulampalam2002} if the effective sample size $M_{Eff}$ in (\ref{eq:EffectiveSampleSize}) is less than $90\%$ of the number of particles. 

\begin{figure}[htbp]
\begin{center}
{\includegraphics[width=1\textwidth]{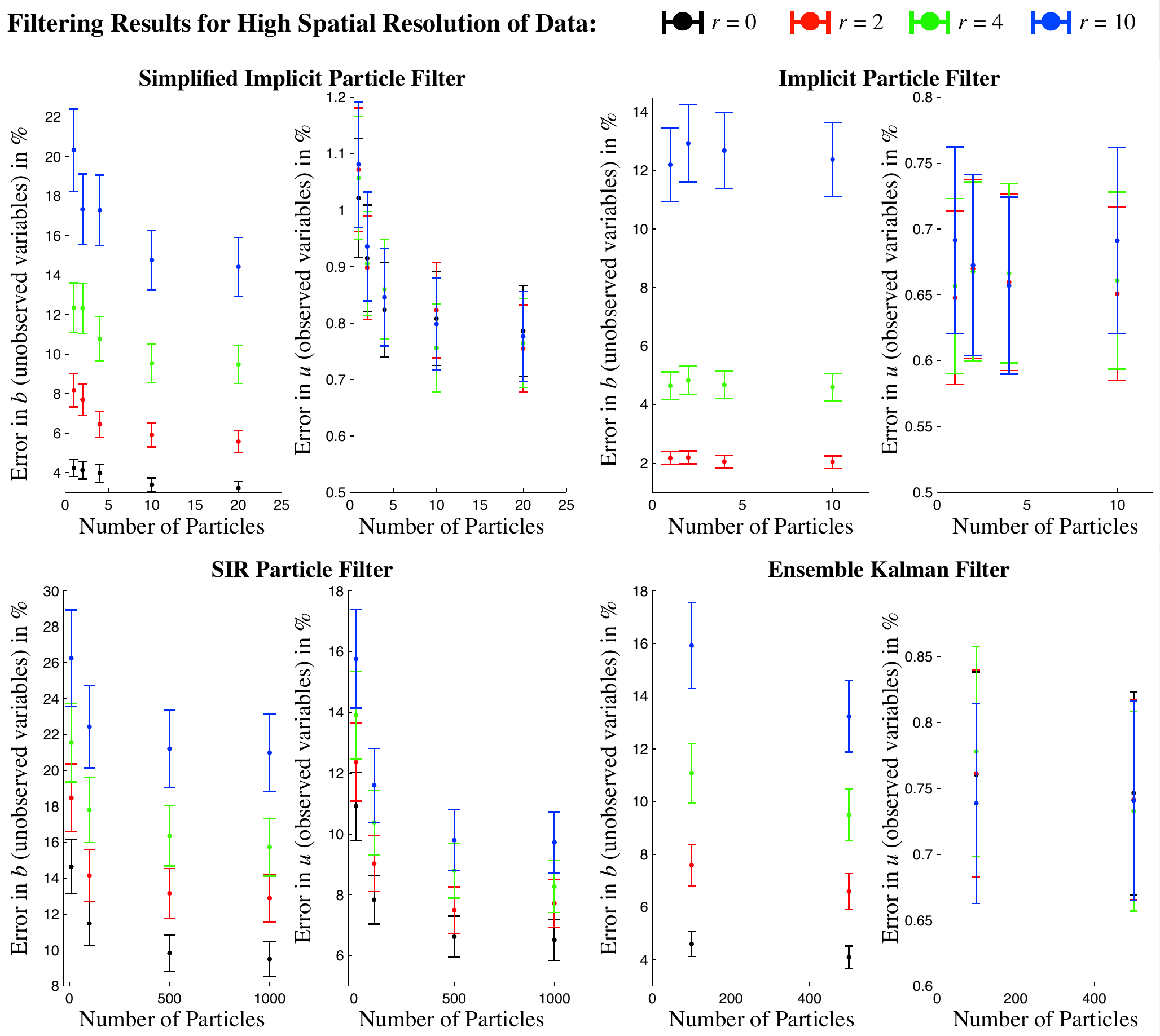}}
\caption{Filtering results for data collected at a high spatial resolution (200 measurement locations). The errors at $T=0.2$ of the simplified implicit particle filter (upper left), implicit particle filter (upper right), SIR filter (lower left) and EnKF (lower right) are plotted as a function of the number of particles and for different gaps between observations in time. The error bars represent the mean of the errors and mean of the standard deviations of the errors.} 
\label{fig:ResultsDense} 
\end{center}
\end{figure}

To further investigate the performance of the filters, we run more numerical experiments and vary the availability of the data in time, as well as the number of particles. Figure \ref{fig:ResultsDense} shows the results for the implicit particle filter, the simplified implicit particle filter, the EnKF and the SIR filter for $200$ measurement locations and for $r=1,2,4,10$.

We observe from Figure \ref{fig:ResultsDense}, that the error statistics of the implicit particle filter have converged, so that there is no significant improvement when we increase the number of particles to more than 10. In fact, the numerical experiments suggest that no more than 4 particles are required here. Independent of the gap between the observations in time, we observe an error of less than 1\% in the observed variable. The error in the unobserved variable $u$ depends strongly on the gap between observations and, for a large gap, is about 12\%. 

The reconstructions of the observed variables by the simplified implicit particle filter are rather insensitive to the availability of data in time and, with 20 particles, the simplified filter gives an error in $u$ of less than 1\%. The errors in the unobserved quantity depend strongly on the gap between the observations and can be as large as 15\%. Here, we need more particles, observe a larger error and a larger sensitivity of the errors to the availability of the data, because the simplified filter makes less direct use of the data, than the ``full'' implicit filter, since it generates the state trajectories using the model and only the final position of each particle is updated by the data. Thus, the error increases as the gap in time between the observations becomes larger. Again, the error statistics have converged and only minor improvements can be expected if the number of particles is increased beyond~20. 

The SIR filter also makes less efficient use of the data so that we require significantly more particles, observe larger errors as well as a stronger dependence of the errors on the availability of data in time, for both the observed and unobserved quantities. With 1000 particles and for a large gap ($r=10$), the SIR filter gives mean errors of 10\% for the observed quantity and 22\% for the unobserved quantity. An increase in the number of particles did not decrease these errors. The EnKF performs well and, for about 500 particles, gives results that are comparable to those of the implicit particle filter. The reason for the large number of particles is, again, the indirect use of the data in EnKF.

\begin{figure}[htbp]
\begin{center}
{\includegraphics[width=1\textwidth]{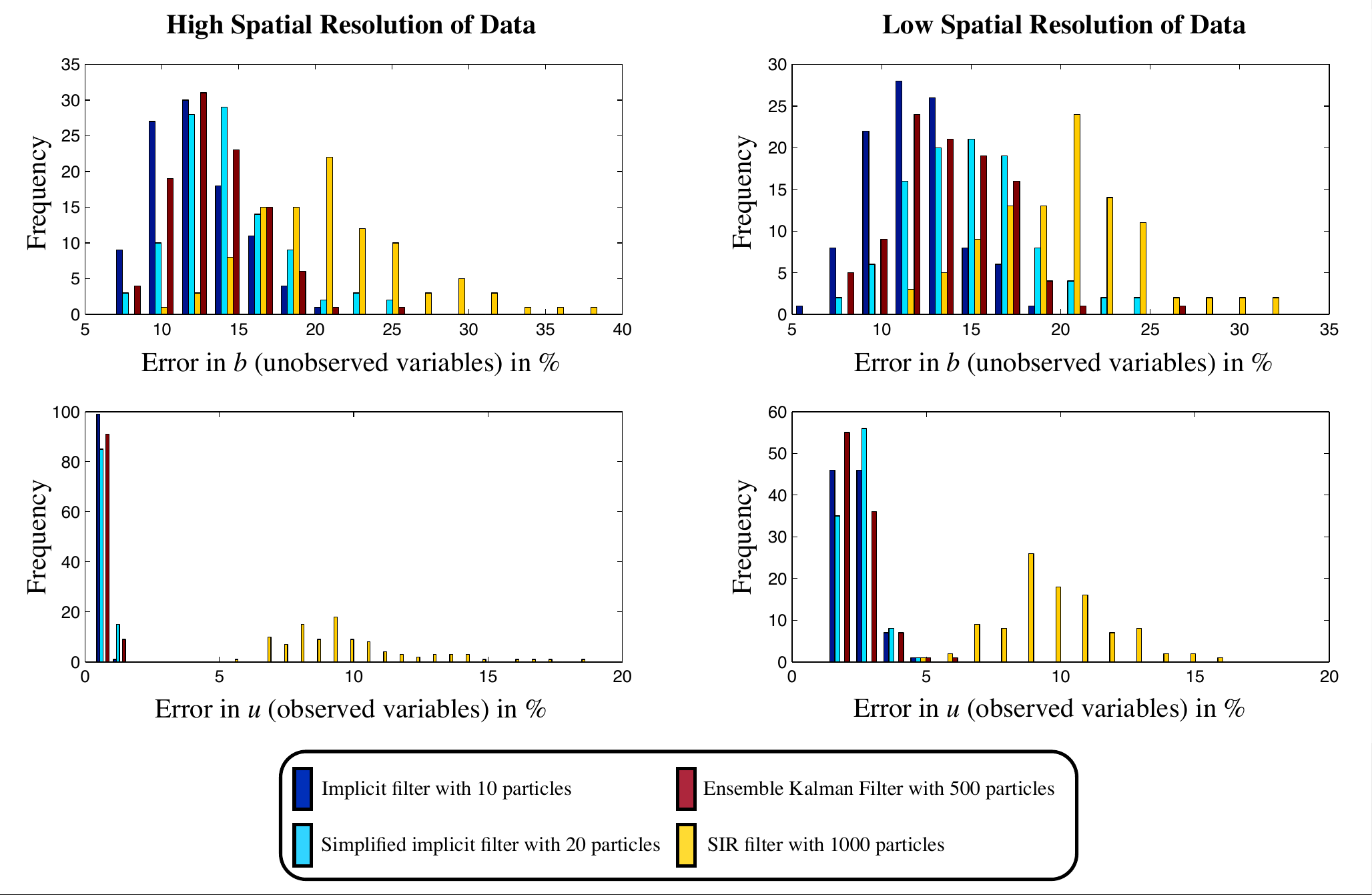}}
\caption{Histogram of errors at $T=0.2$ of the implicit filter, simplified implicit filter, EnKF and SIR filter. Left: data are available at a high spatial resolution (200 measurement locations) and every $r=10$ model steps. Right: data are available at a low spatial resolution (20 measurement locations) and every $r=10$ model steps.} 
\label{fig:Distribution} 
\end{center}
\end{figure}

The errors in the reconstructions of the various filters are not Gaussian, so that an assessment of the errors based on the first two moments is incomplete. In the two panels on the right of Figure \ref{fig:Distribution}, we show histograms of the errors of the implicit filter (10 particles), simplified implicit filter (20 particles), EnKF (500 particles) and SIR filter (1000 particles) for $r=10$ model steps between observations. 
We observe that the errors of the implicit filter, simplified implicit filter and EnKF are centered to the right of the diagrams (at around 10\% in the unobserved quantity $u$ and about 1\% for the observed quantity $b$) and show a considerably smaller spread than the errors of the SIR filter, which are centered at much larger errors (20\% in the unobserved quantity $u$ and about 9\% for the observed quantity $b$). A closer look at the distribution of the errors thus confirms our conclusions we drew from an analysis based on the first two moments. 

\begin{figure}[htbp]
\begin{center}
{\includegraphics[width=1\textwidth]{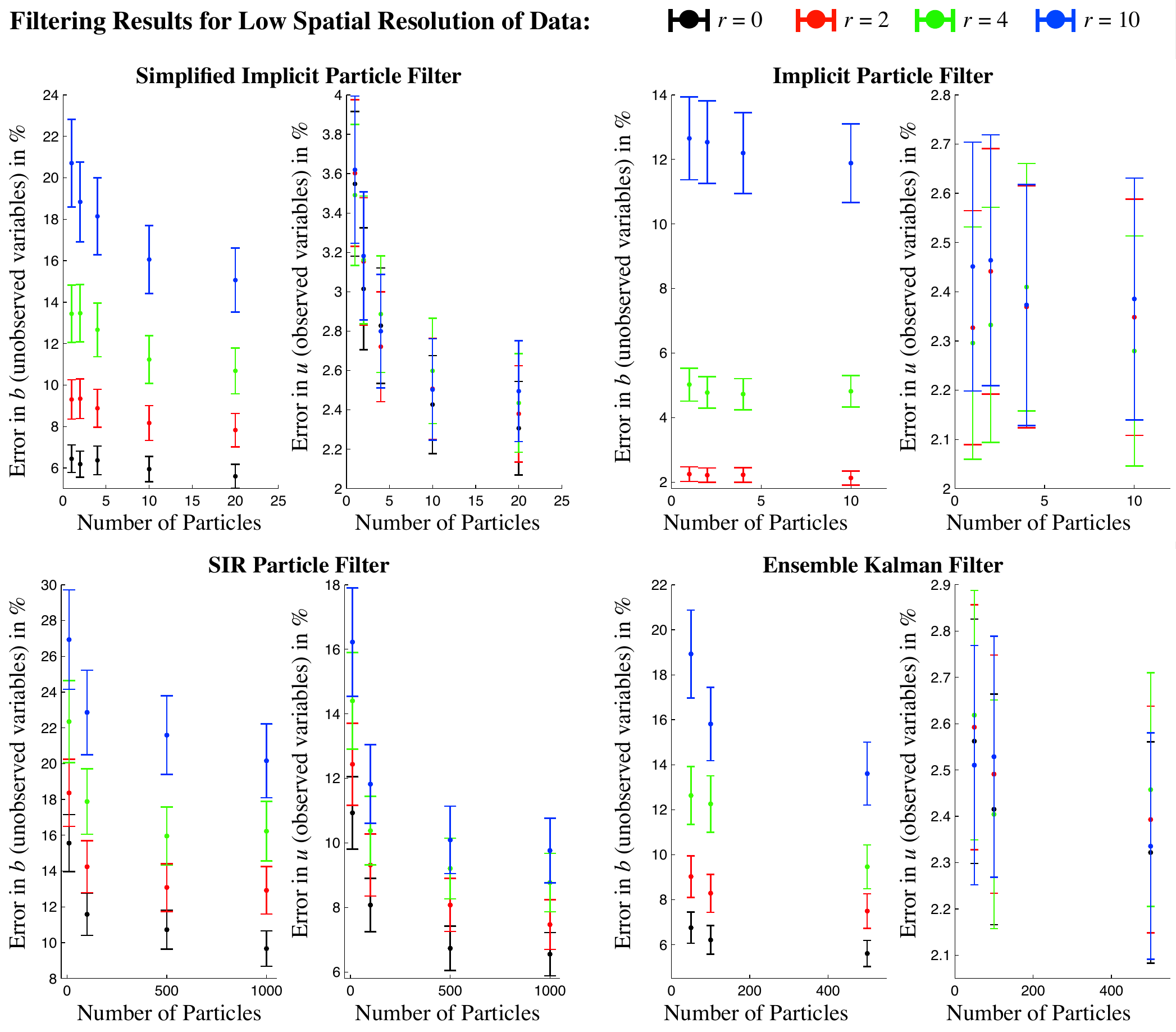}}
\caption{Filtering results for data collected at a low spatial resolution (20 measurement locations). The errors at $T=0.2$ of the simplified implicit particle filter (upper left), implicit particle filter (upper right), SIR filter (lower left) and EnKF (lower right) are plotted as a function of the number of particles and for different gaps between observations in time. The error bars represent the mean of the errors and mean of the standard deviations of the errors.} 
\label{fig:ResultsSparse} 
\end{center}
\end{figure}

We decrease the spatial resolution of the data to $20$ measurement locations and show filtering results from 100 twin experiments in Figure \ref{fig:ResultsSparse}.
The results are qualitatively similar to those obtained at a high spatial resolution of 200 data points per observation. We observe for the implicit particle filter that the errors in the unobserved quantity are insensitive to the spatial resolution of the data, while the errors in the observed quantity are determined by the spatial resolution of the data and are rather insensitive to the temporal resolution of the data. These observations are in line with those reported in connection with a strong 4D-Var algorithm in \cite{Fournier2007}. All other filters we have tried show a dependence of the errors in the observed quantity on the temporal resolution of the data. Again, the reason for the good performance of the implicit particle filter is its direct use of the data. The two panels to the left of Figure \ref{fig:Distribution}, show histograms of the errors of the implicit filter (10 particles), simplified implicit filter (20 particles), EnKF (500 particles) and SIR filter (1000 particles) for $r=10$ model steps between observations. The results are qualitatively similar to the results we obtained at a higher spatial resolution of the data and the closer look at the distributions of the errors confirms the conclusions we drew from an analysis of the first two moments.

In summary, we observe that the implicit particle filter yields the lowest errors with a small number of particles for all examples we considered, and performs well and reliably in this application. The SIR and simplified implicit particle filters could not reach the accuracy of the implicit particle filter, even when the number of particles is very large. The EnKF requires about 500 particles to come close to the accuracy of the implicit particle filter with only 4 particles. Although the implicit filter uses the computationally most expensive particles, the small number of particles required for a very high accuracy make the implicit filter the most efficient filter for this problem. The partial noise works in our favor here, because the dimension of the space the implicit filter operates in is 20, rather than the state dimension 600. 

Finally, we wish to compare our results with those in \cite{Fournier2007}, where a strong constraint 4D-Var algorithm was applied to the deterministic version of the test problem. Fournier used ``perfect data,'' i.e. the observations were not corrupted by noise, and applied a conjugate-gradient algorithm to minimize the 4D-Var cost function. The iterative minimization was stopped after 5000 iterations. With 20 observations in space and a gap of $r=5$ model steps between observations, an error of about 1.2\% in $u$ and 4.7\% in $b$ was achieved. With the implicit filter, we can get to a similar accuracy at the same spatial resolution of the data, but with a larger gap of $r=10$ model steps between observations. Moreover, the data assimilation problem we solve here is somewhat harder than the strong constraint 4D-Var problem because we allow for model errors. The implicit particle filter also reduces the memory requirements because it operates in the 20-dimensional subspace of the forced variables. Each minimization is thus not as costly as a 600-dimensional strong constraint 4D-Var minimization.

\section{Conclusions}
We considered implicit particle filters for data assimilation. Previous implementations of the implicit particle filter rely on finding the Hessians of functions $F_j$ of the state variables. Finding these Hessians can be expensive if the state dimension is large and can be cumbersome if the second derivatives of the $F_j$'s are hard to calculate. We presented a new implementation of the implicit filter combining gradient descent minimization with random maps. This new implementation avoids the often costly calculation of the Hessians and, thus, reduces the memory requirements compared to earlier implementations of the filter. 

We have considered models for which the state covariance matrix is singular or ill-conditioned. This happens often, for example, in geophysical applications in which the noise is smooth in space or if the model includes conservation laws with zero uncertainty. Previous implementations of the implicit filter are not applicable here and we have shown how to use our new implementation in this situation. The implicit filter is found to be more efficient than competing methods because it operates in a space whose dimension is given by the rank of the state covariance matrix rather than the model dimension.
 
We applied the implicit filter in its new implementation to a test problem in geomagnetic data assimilation. The implicit filter performed well in comparison to other data assimilation methods (SIR, EnKF and 4D-Var) and gave accurate state estimates with a small number of particles and at a low computational cost. We have studied how the various data assimilation techniques use the available data to propagate information from observed to unobserved quantities and found that the implicit particle filter uses the data in a direct way, propagating information to unobserved quantities faster than competing methods. The direct use of the data is the reason for the very small errors in reconstructions of the state.

\section*{Acknowledgments}
We would like to thank our collaborators Dr. Ethan Atkins at UC Berkeley, and Professors Robert Miller, Yvette Spitz, and Dr. Brad Weir at Oregon State University, for their comments and helpful discussion. We thank Mr. Robert Saye for careful proofreading of early versions of this manuscript. This work was supported in part by the Director, Office of Science, Computational and Technology Research, U.S. Department of Energy under Contract No. DE-AC02-05CH11231.

\bibliographystyle{plain}
\bibliography{References}

\end{document}